\documentclass[10pt,a4paper]{amsart}
\usepackage[shortlabels]{enumitem}
\usepackage{amsmath,amssymb,amsxtra,amsthm,stackrel,mathrsfs}
\usepackage{mathpazo}
\usepackage[normalem,normalbf]{ulem}
\usepackage[all]{xy}
\usepackage[british]{babel}
\usepackage[utf8]{inputenc}
\usepackage[colorlinks=true,linkcolor=blue,citecolor=blue]{hyperref}

\newtheorem{prp}{Proposition}[section]
\newcommand{\bp}{\begin{prp}}
\newcommand{\ep}{\end{prp}}
\newtheorem{lemme}[prp]{Lemma}
\newcommand{\bl}{\begin{lemme}}
\newcommand{\el}{\end{lemme}}
\newtheorem{cor}[prp]{Corollary}
\newcommand{\bc}{\begin{cor}}
\newcommand{\ec}{\end{cor}}
\newtheorem{thm}[prp]{Theorem}
\newcommand{\bt}{\begin{thm}}
\newcommand{\et}{\end{thm}}

\theoremstyle{remark}
\newtheorem*{rmk}{Remark}
\newcommand{\br}{\begin{rmk}}
\newcommand{\er}{\end{rmk}}

\newcommand{\bpf}{\begin{proof}}
\newcommand{\epf}{\end{proof}}

\newcommand{\ppq}{\leqslant}
\newcommand{\pgq}{\geqslant}
\newcommand{\ot}{\otimes}
\newcommand{\B}{\mathcal{B}}
\newcommand{\A}{\mathcal{A}}
\newcommand{\C}{\mathcal{C}}
\renewcommand{\L}{\mathcal{L}}
\newcommand{\D}{\mathcal{D}}
\newcommand{\id}{\operatorname{id}\nolimits}
\newcommand{\im}{\operatorname{Im}\nolimits}
\renewcommand{\ker}{\operatorname{Ker}\nolimits}
\newcommand{\trc}{\operatorname{trace}\nolimits}
\newcommand{\der}{\operatorname{\mathsf{der}}\nolimits}
\newcommand{\Hom}{\operatorname{Hom}\nolimits}
\newcommand{\Ext}{\operatorname{Ext}\nolimits}
\newcommand{\HH}{\operatorname{HH}\nolimits}
\newcommand{\sHH}{\operatorname{\uline{HH}}\nolimits}
\newcommand{\sEnd}{\operatorname{\uline{End}}\nolimits}
\newcommand{\End}{\operatorname{End}\nolimits}
\newcommand{\car}{\operatorname{char}\nolimits}
\newcommand{\mg}{\ensuremath{\mathfrak{g}}}
\newcommand{\rad}{\operatorname{rad}\nolimits}
\newcommand{\s}[1]{\stackrel[#1]{}{\otimes}}
\newcommand{\set}[1]{\left\{ #1 \right\}}
\newcommand{\rep}[1]{\langle{#1}\rangle} 
\newcommand{\pbar}{\operatorname{Bar}(A)} 
\newcommand{\pmin}{\mathcal{P}}
\newcommand{\ov}{\overline}
\newcommand{\tup}{\textsuperscript}
\newcommand{\spn}[1]{\operatorname{span}\set{#1}}

\title[HH1 and stable equivalence of Morita type]{First Hochschild cohomology group and stable equivalence classification of Morita type of some tame symmetric algebras}
\author[R. Taillefer]{Rachel Taillefer}
\address{R. Taillefer, Universit\'e Clermont Auvergne, CNRS, LMBP,  F-63000 Clermont-Ferrand, France}
\email{Rachel.Taillefer@uca.fr}
\date{\today}
\subjclass[2010]{
17B60,   %Lie (super)algebras associated with other structures (associative, Jordan, etc.)
16E40, % (Co)homology of rings and algebras (e.g. Hochschild, cyclic, dihedral, etc.)
16G99, % None of the above, but in this section: Representation theory of rings and algebras
16D90.   %Module categories [See also 16Gxx, 16S90]; module theory in a category-theoretic context; Morita equivalence and duality
}

\keywords{stable equivalence of Morita type; algebras of dihedral, semi-dihedral and quaternion type; first Hochschild cohomology group; Lie algebra.}

\begin{document}

\maketitle
\begin{abstract}
 We use the dimension and the Lie algebra structure of the first Hochschild cohomology group to distinguish some algebras of dihedral, semi-dihedral and quaternion type up to stable equivalence of Morita type. In particular, we complete the classification of algebras of dihedral type that was mostly determined by Zhou and Zimmermann.
\end{abstract}

\section*{Introduction}
Erdmann has given a description, up to Morita equivalence, of some families of tame symmetric
algebras, which include the blocks of finite group algebras of tame representation type, and that are defined essentially in terms of their Auslander-Reiten quivers. They are separated into three types, dihedral, quaternion and semi-dihedral (generalising tame blocks whose defect groups are dihedral, semi-dihedral or generalised quaternion).  Holm then classified them up to derived equivalence in \cite{Ho}. It is then natural to try to classify them up to stable equivalence, but there are many properties that are not preserved under stable equivalences. However, Rickard in \cite{Ri} and Keller and Vossieck in \cite{KV} proved that a derived equivalence between selfinjective algebras induces a stable equivalence of a particular form, called stable equivalence of Morita type because it is induced by tensoring with some bimodules; since then, such stable equivalences (even for algebras that are not selfinjective) have been much studied. In particular, 
in \cite{ZZ} and in \cite{Z},  Zhou and  Zimmermann used various techniques (including K\"ulshammer invariants and
stable Hochschild cohomology)
in order to distinguish most of the algebras of dihedral, semi-dihedral and quaternion type up to stable equivalence of Morita type, but some questions remain. Our aim is to use  the first Hochschild
cohomology group and its Lie structure to answer some of these questions.

 It was shown by Xi in \cite{X} that if $A$ and $B$ are two selfinjective algebras
and if there is a stable equivalence of
Morita type between them, then 
for $n\pgq1$,  the
Hochschild cohomology groups $\HH^n(A)$ and $\HH^n(B)$ are isomorphic. Moreover, as a consequence of
a result of 
 K\"onig,  Le and  Zhou  in \cite{KLZ},  if $A$ is a symmetric algebra, the Lie algebra structure of $\HH^1(A)$ is also
preserved under such an equivalence. We shall use these facts to distinguish some of the algebras
above up to stable equivalence
of Morita type. As a result, we are able to complete the classification for the  algebras of dihedral
type, and to improve it for the  algebras of quaternion and semi-dihedral types.  

The Lie algebra structure of the first Hochschild cohomology group has been described by Strametz in \cite{S}, where she studied the Lie algebra $\HH^1(A)$  for a monomial algebra $A$. Her results were then used by Sánchez-Flores in \cite{SF} to study the Gerstenhaber algebra structure of the Hochschild cohomology ring $\HH^*(A)$ of a monomial algebra $A$.  Strametz' description has also been used by Bessenrodt and Holm in \cite{BH}. The Lie algebra $\HH^1(A)$ has also been studied for instance in \cite{GAS}, and used for example in \cite{LRD} to retrieve information on some blocks of a group algebra. We shall describe Strametz' construction in Section \ref{sec:general description Lie algebra} and use it in this paper.

After summarising in Section \ref{sec:questions} the results known on stable equivalence of Morita type of algebras of dihedral, semi-dihedral and quaternion type, as well as proving our main result for algebras of quaternion type with two simple modules, we give some general tools that we will use in Section \ref{sec:general}: we first describe the Lie algebra structure on the first Hochschild cohomology group. Moreover, the usual algorithmic methods to compute a minimal projective resolution of an algebra given by quiver and relations relies on the fact that we have a minimal set of relations, which is not the case here. Therefore we describe our method to determine the beginning of a minimal projective resolution of a finite-dimensional associative algebra in order to compute the first Hochschild cohomology group. Finally, we shall use some constructions that are invariant under Lie algebra isomorphisms, which we recall in the last part of  Section \ref{sec:general}. We then study the cases of algebras of dihedral type in Section \ref{sec:dihedral},   of semi-dihedral type in Section \ref{sec:semidihedral},  and of quaternion type in Section \ref{sec:quaternion}.

 Throughout, $K$ is an algebraically closed field. Set $\ot=\ot_K.$

\subsection*{Acknowledgements.}\sloppy I wish to thank Alexander Zimmermann for asking me, many years ago, whether the dimension of $\HH^1(A)$ could help with the classification problem for the local dihedral algebras, thus initiating this project.

\section{The questions studied in this paper}\label{sec:questions}

In \cite{ZZ}, Zhou and Zimmermann proved that if $A$ and $B$ are algebras that are stably equivalent of Morita type, then $A$ is of dihedral (respectively semi-dihedral, respectively quaternion) type if and only if $B$ is also. Moreover, if $A$ and $B$ are of dihedral, semi-dihedral and quaternion type, then $A$ and $B$ have the same number of simple modules. 

Since our methods did not enable us to improve on the existing results for algebras with three simple modules (the only question being for the algebras $Q(3\A)_1^{2,2}(d)$  of quaternion type for which the Lie algebra structure of $\HH^1(Q(3\A)_1^{2,2}(d))$ does not depend on $d$), we shall restrict our study to the algebras with one or two simple modules. 

\subsection{The algebras involved}\label{subsec:questions:list algebras}

Let us first define the $K$-algebras that we are going to study, by quiver and relations. We shall need the following quivers:
\[ \xymatrix{
\bullet\ar@(dl,ul)[]^x\ar@(dr,ur)[]_y&&&&&1\ar@/^/[rr]^\beta\ar@(dl,ul)[]^\alpha&&2\ar@/^/[ll]^\gamma\ar@(dr,ur)[]_\eta\\
1\A  &&&&&&2\B 
} \] 

We shall only consider the local algebras when $\car(K)=2$, and they are defined as follows. 

The quiver of all the local algebras is $1\A$. Moreover,
\begin{itemize}
\item the algebras $D(1\A)_2^k(d)$ of dihedral type, where $k\pgq2$ is an integer and $d\in\set{0,1}$, whose relations ideal is generated by \[x ^2 - (xy )^k ,\  y ^2 - d  (xy )^k ,\  (xy )^k - (yx )^k ,\  (xy )^k x \text{ and } (yx )^ky,\]
\item  the algebras $SD(1\A)_2^k(c,d)$ of semi-dihedral type, where $k\pgq2$ is an integer and $(c,d)\in K^2$, $(c,d)\neq (0,0)$,  whose relations ideal is generated by \[(xy )^k - (yx )^k ,\  (xy )^k x ,\  y ^2 - d(xy )^k\text{ and } x ^2 - (yx )^{k-1}y + c(xy )^k ,\]
\item  the algebras $Q(1\A)_1^k$ of quaternion type, where $k\pgq2$ is an integer,  whose relations ideal is generated by  \[(xy )^k - (yx )^k ,\  (xy )^kx ,\  y ^2 - (xy )^{k-1} x \text{ and } x ^2 - (yx )^{k-1} y \] and the algebras $Q(1\A)_2^k(c,d)$ of quaternion type, where $k\pgq2$ is an integer and $(c,d)\in K^2$, $(c,d)\neq (0,0)$, whose relations ideal is generated by \[x ^2 - (yx )^{k-1} y - c (xy )^k,\  y ^2 - (xy )^{k-1} x - d(xy )^k ,\  (xy )^k - (yx )^k ,\  (xy )^k x  \text{ and }(yx )^k y.\]
\end{itemize}

These algebras all have dimension $4k$ with basis  the elements 
\begin{align*}
&(xy)^t&&(yx)^{t+1}&&y(xy)^t&&x(yx)^t
\end{align*}
for $0\ppq t\ppq k-1$, and the centre of all these algebras has dimension $(k+3)$.

We no longer assume that $\car(K)=2$. The quiver of all the algebras with two simple modules is $2\B$ and they are the following.
\begin{itemize}
\item  The algebras $SD(2\B)_1^{k,s}(c)$ of semi-dihedral type, where $k\pgq2$ and $s\pgq 1$ are integers and $c\in\set{0,1}$,  whose relations ideal is generated by \[\gamma  \beta ,\  \eta \gamma  ,\  \beta \eta , \alpha  ^2 - (\beta \gamma  \alpha )^{k-1} \beta \gamma  - c (\alpha \beta \gamma  )^k ,\  \eta^ s - (\gamma  \alpha \beta )^k\text{ and } (\alpha \beta \gamma  )^k - (\beta \gamma  \alpha )^k.\]
\item  The algebras $SD(2\B)_2^{k,s}(c)$ of semi-dihedral type, where $k\pgq2$ and $s\pgq 1$ are integers  with $k+s\pgq 4$ and $c\in\set{0,1}$,  whose relations ideal is generated by \[\beta \eta  - (\alpha \beta \gamma  )^{k-1} \alpha \beta ,\  \eta \gamma  - (\gamma  \alpha \beta )^{k-1} \gamma  \alpha ,\  \gamma  \beta  - \eta ^{s-1},\ \alpha  ^2 - c (\alpha \beta \gamma  )^k ,\  \beta \eta ^2\text{ and } \eta ^2 \gamma  .\]
\item  The algebras $Q(2\B)_1^{k,s}(a,c)$ of quaternion type, where $k\pgq1$ and $s\pgq 3$ are integers and $(a,c)\in K^2$ with $a\neq 0$ (and $a\neq 1$ if $k+s=4$),  whose relations ideal is generated by 
\begin{align*}
&\gamma  \beta  - \eta  ^{s-1},&& \beta \eta  - (\alpha \beta \gamma  )^{k-1} \alpha \beta ,&&  \eta \gamma  - (\gamma  \alpha \beta )^{k-1} \gamma  \alpha ,\\&  \alpha ^2 - a(\beta \gamma  \alpha )^{k-1}\beta \gamma  - c(\beta \gamma  \alpha )^k,&&  \alpha  ^2\beta,&&\gamma\alpha^2.
\end{align*}
\end{itemize}

The algebras with two simple modules (of semi-dihedral and quaternion type) all have  dimension $9k+s$, and  the following elements, where $0\ppq t\ppq k-1$ and $1\ppq r\ppq s$, form a basis of each algebra:
\begin{xalignat*}{5}
&(\alpha\beta\gamma)^t, && (\beta\gamma\alpha)^{t+1},&&(\alpha\beta\gamma)^t\alpha, &&
(\beta\gamma\alpha)^t\beta\gamma,&&(\beta\gamma\alpha)^t\beta,\\&(\alpha\beta\gamma)^t\alpha\beta,&&
(\gamma\alpha\beta)^t\gamma,&&(\gamma\alpha\beta)^t\gamma\alpha,&&(\gamma\alpha\beta)^{t+1}, &&\eta^r.
\end{xalignat*} 
Moreover, their centre has dimension $k+s+2.$

\subsection{Algebras of dihedral type}\label{subsec:questions:dihedral}

In the case of algebras of dihedral type, Zhou and Zimmermann proved that the classification up to stable equivalence of Morita type mostly coincides with the classification up to derived equivalence, but a few questions in the classification remain. As they stated in \cite[Remark 4.2 and Remark 7.2]{ZZ},  in order to complete the classification of the algebras of dihedral type we must determine whether the algebras $D(1\A)_2^k(0)$ and  $D(1\A)_2^k(1)$ are stably equivalent of Morita type or not. We shall prove that they are not in Corollary \ref{cor:local dihedral}.

\subsection{Algebras of semi-dihedral type}\label{subsec:questions:semidihedral}

\sloppy The remaining question for the local algebras of semi-dihedral type is whether the  stable equivalence of Morita type classes for the algebras $SD(1\A)_2^k(c,d)$ depend on $(c,d)$ or not. We shall give a partial answer to this question in  Theorem \ref{thm:classification sd local}.

\medskip

In the case of algebras of semi-dihedral type with two simple modules, it was proved in \cite{ZZ} that if two such algebras $SD(2\B)_i^{k,s}(c)$ and $SD(2\B)_{i'}^{k',s'}(c')$ are stably equivalent of Morita type, then the sets $\set{k,s}$ and $\set{k',s'}$ are equal. Moreover, if $\car(K)=2$ and  $k=2$ and $s\pgq 3$ is odd then  $SD(2\B)_1^{k,s}(0)$ and $SD(2\B)_{1}^{k',s'}(1)$ are not stably equivalent of Morita type, and  if $k$ and $s$ are both odd, and if $\set{k',s'}=\set{k,s}$, then  $SD(2\B)_2^{k,s}(0)$ and $SD(2\B)_{2}^{k',s'}(1)$ are not stably equivalent of Morita type.

We go further in this classification in Theorem \ref{thm:sd}.

\subsection{Algebras of quaternion type}\label{subsec:questions:quaternion}

For local algebras of quaternion type, the remaining questions are whether, for a fixed $k$, two algebras among $Q(1\A)_1^k$ and the $Q(1\A)_2^k(c,d)$ can be stably equivalent of Morita type or not. We shall study this situation in Section \ref{sec:quaternion}, whose main result is Corollary \ref{cor:classification q 2}.

\medskip

We now turn to the algebras of quaternion type with two simple modules.

 Zhou and Zimmermann showed in \cite{ZZ} that if  $Q(2\B)_1^{k,s}(a,c)$ and  $Q(2\B)_1^{k',s'}(a',c')$ are stably equivalent, then the sets $\set{k,s}$ and $\set{k',s'}$ are equal. There remained some questions however.

First assume that  $\car(K)\neq 2$. If $k+s>4$, the problem relating to the parameters $a$ and $c$ was solved recently by Zimmermann in \cite{Z}, where he proved that  $Q(2\B)_1^{k,s}(a,c)\cong  Q(2\B)_1^{k,s}(1,0) $ (recall that the field $K$ is algebraically closed). Therefore, using \cite[Lemma 5.7 (ii)]{ES},  we need only consider the algebras $Q(2\B)_1^{1,3}(a,0)$ with $a\in K^*$, $a\neq 1.$ However, the methods in this paper do not provide any information to distinguish the stable equivalence classes of Morita type, therefore we shall assume that $\car(K)=2.$
 In this case,  then by \cite[Lemma 5.7]{ES}, if $k+s>4$  we need only consider the algebras $Q(2\B)_1^{k,s}(1,c)$ for $c\in K$. 

Theorem \ref{thm:quaternion 2 simples} below can be obtained from the Lie algebra structure of $\HH^1(Q(2\B)_1^{k,s}(a,c))$ and the techniques in this paper, using a minimal projective resolution from \cite{ES} and computing the dimensions of the Hochschild cohomology groups as in Proposition \ref{prop:HHn quaternion 1 simple} and the Lie algebra structure of the first Hochschild cohomology group as in the other cases. However we shall give a more elegant proof here using a result from \cite{Z}. 

\bt\label{thm:quaternion 2 simples} Let $K$ be an algebraically closed field of characteristic $2$.  Let $k\pgq 1$ and $s\pgq3$ be integers,  and let $c$ be in $K^*.$ 

 If $k+s>4$, then the algebras    $Q(2\B)_1^{k,s}(1,c)$ and $Q(2\B)_1^{k,s}(1,0)$ are not stably equivalent of Morita type. 

If $k=1$ and $s=3$  then, for any $a,a'$ in $K\backslash\set{0,1}$, the algebras    $Q(2\B)_1^{k,s}(a,c)$ and $Q(2\B)_1^{k,s}(a',0)$ are not stably equivalent of Morita type. 
\et

Before we prove this result, let us define the objects that we shall use. Let $A$ be a symmetric algebra over a field of characteristic $p$, endowed with a non-degenerate symmetric associative bilinear form $(,)$. Let $KA$ be the subspace of $A$ generated by the commutators $ab-ba$ of elements $a$, $b$ in $A$. Define $T_n(A)=\set{x\in A,\ x^p\in KA}$ and let $T_n(A)^\perp$ be the orthogonal space with respect to $(,)$, which is an ideal in the centre $Z(A)$, called Külshammer ideal. 

The algebra $Z(A)/T_1^\perp(A)$ is a stable invariant of Morita type. Indeed, 
 let
  $Z^{st}(A)=\sEnd_{A^e}(A)$ be the
  stable centre of $A$ (the endomorphisms of $A$ in the stable $A^e$-module category)  and let  $Z^{pr}(A)=\ker({\End}_{A^e}(A)\rightarrow\sEnd_{A^e}(A))$ be
  the projective centre of $A.$ Then the ideals $Z^{st}(A)$ and $T_1^\perp(A)/Z^{pr}(A)$ are invariants of
  stable equivalences of Morita type for symmetric algebras (see \cite{LZZ,KLZ}), and  moreover
  $Z(A)/T_1^\perp(A)\cong Z^{st}(A)/(T_1^\perp(A)/Z^{pr}(A))$.

\bpf 
In \cite[Theorem 7 (2)]{Z}, Zimmermann describes the quotient  $Z(Q(2\B)_1^{k,s}(a,c))/T_1^\perp(Q(2\B)_1^{k,s}(a,c))$ in all cases depending on the values and parity of $k$ and $s$ and on whether $c=0$ or $c\neq 0$,  and  it follows that the algebras
 $Z(Q(2\B)_1^{k,s}(a,c))/T_1^\perp(Q(2\B)_1^{k,s}(a,c))$  and  $Z(Q(2\B)_1^{k,s}(a',0))/T_1^\perp(Q(2\B)_1^{k,s}(a',0))$ are not isomorphic when $c\neq 0$ (see the proof of \cite[Corollary 10]{Z}). The result follows.
\epf

\br The same result when $k\pgq 2$ can be obtained as a consequence of the algebra structure of the whole Hochschild cohomology computed in \cite{GII}. We note that although the algebras $\HH^*(Q(2\B)_1^{k,s}(1,c))$ and $\HH^*(Q(2\B)_1^{s,k}(1,c))$ in \cite{GII} appear to be different, there is an explicit isomorphism between them.
\er

\br 
We should  mention that the first Hochschild cohomology group does not separate algebras with different parameter $a.$
\er

Thereofore, in Section \ref{sec:quaternion} we shall only consider local algebras of quaternion type.

\section{General facts on the first Hochschild cohomology group and its computation and on invariants of Lie algebras}\label{sec:general}

\subsection{Lie algebra structure on $\HH^1(A)$}\label{sec:general description Lie algebra}

\sloppy K\"onig, Le and  Zhou proved in \cite[Theorem 10.7]{KLZ} that the Batalin-Vilkoviskyi
structure of the stable Hochschild
cohomology $\HH_{st}^*(A)$ (that is, the Hochschild cohomology $\HH^*(A)=\bigoplus_{n\in\mathbb{N}}\HH^n(A)$
modulo the projective centre of $A$) of a symmetric algebra $A$ is invariant under stable equivalences of
Morita type. In particular, the Lie algebra structure of $\HH^1(A)$ is preserved under such an equivalence. 

The Lie structure on $\HH^1(A)$ is usually described on the Hochschild complex (obtained from the bar resolution). However, we will be working with minimal resolutions, so we will need a description of the Lie bracket when  $\HH^1(A)$ is computed from a minimal projective resolution. This is based on \cite{S}.

Let $A=K\Gamma/I$ be a finite dimensional algebra, where $\Gamma$ is a quiver and $I$ is an admissible
ideal.
Let $\Gamma_0$ be the set of vertices in $\Gamma$ and $\Gamma_1$ be the set of arrows.

Using the bar resolution $\pbar$, we get $\HH^{1}(A)=\ker d_1/\im d_0$ where 
\[ 0\rightarrow  A\xrightarrow {d_0}{} \Hom_K(A,A)\xrightarrow {d_1}{}\Hom_K(A\ot _KA,A)  \]
and $d_0(\lambda)(p)=\lambda p-p\lambda$ and $d_1(f)(a\ot b)=af(b)-f(ab)+f(a)b.$

The space $\HH^1(A)$ is then endowed with a Lie bracket defined by 
\[ [f,g]=f\circ g-g\circ f. \]

For all the algebras $A$ we shall consider in this paper, there is a minimal projective resolution $\pmin$ of $A$ that starts with 
\[ A\ot _E KZ\ot _EA\stackrel{\partial^1}{\longrightarrow} A\ot _E K\Gamma _1\ot
_EA\stackrel{\partial^0}{\longrightarrow} A\ot _E A\rightarrow A\rightarrow 0 \]
where $E=K\Gamma _0$, $Z$ is a set of relations in $I$ and 
\begin{align*}
\partial^0(1\ot_Ea\ot_E1)&=a\ot_E1-1\ot _Ea\qquad\text{ for $a\in\Gamma_1$}\\
\partial^1(1\ot_E\left(\sum_{i=1}^nc_ia_{1,i}\cdots
  a_{s_i,i}\right)\ot_E1)&=\sum_{i=1}^n\sum_{j=1}^{s_i}c_ia_{1,i}\cdots
a_{j-1,i}\ot_Ea_{j,i}\ot_Ea_{j+1,i}\cdots a_{i,s_i}
\end{align*} where the $c_i$ are in $ K$, the $a_{j,i}$ are in $\Gamma_1$ and $z=\sum_{j=1}^{s_i}c_ia_{j,i}\in Z.$

 As $\pbar$ and $\pmin$ are projective 
 resolutions of the $A$-bimodule $A$, there exist, by the 
 Comparison Theorem, chain maps $\xi\colon \pbar\rightarrow
 \pmin$ and 
 $\varrho\colon \pmin\rightarrow
 \pbar$. As in \cite{S}, these maps induce inverse linear isomorphisms $\xi^*_1$ and $\varrho_1^*$  at the cohomology level between $\HH^1(A,A)=
 \ker d_1/\im d_0$ and
 $\HH^1(A,A)=\ker\partial_1/\im\partial_0$  given by the classes of 
 \[\begin{array}{rrcl}
   \ov{\xi_1}\colon &\Hom_{E-E}(K\Gamma _1,A)&\longrightarrow&\Hom_K(A,A)\\
   &f&\longmapsto&
    \left[
     a_1\cdots a_n\mapsto\sum_{i=1}^na_1\cdots a_{i-1}f(a_i)a_{i+1}\cdots a_n\right]
   \\
   \ov{\varrho_1}\colon &\Hom_K(A,A)&\longrightarrow&\Hom_{E-E}(K\Gamma _1,A)\\
    & h&\mapsto&h_{|K\Gamma _1}.
 \end{array}\]
 This allows us to transfer the Lie algebra structure of $\ker d_1/\im d_0$ to
 $\ker\partial_1/\im\partial_0$, where the bracket is given by 
 \[[f,g]\colon ={\varrho_1}^*([\xi_1^*(f),
   \xi_1^*(g)])=\xi_1^*(f)\circ g- \xi_1^*(g)\circ f\] for all $f,g$ in
   $\Hom_{E-E}(K\Gamma _1,A)$.

\subsection{Method used to determine the beginning of a minimal projective resolution of an algebra $A$ as an $A$-$A$-bimodule}

Given a finite-dimensional $K$-algebra $A=K\Gamma/I$ defined by quiver $\Gamma$ and relations $I$, Happel's theorem \cite{H} gives the modules in a  minimal projective resolution of an algebra $A$ as an $A$-$A$-bimodule, but not the maps.  The general methods to determine the  beginning of a minimal projective resolution of an algebra $A$ as an $A$-$A$-bimodule usually rely on the fact we have a minimal set of generators for the algebra $I$. However,  most of the algebras of dihedral, semi-dihedral and quaternion type are not defined with a minimal set of relations, and it is not easy to extract such a minimal set. Therefore, we shall repeatedly use the following result of  \cite[Proposition 2.8]{GS} (see also \cite[Theorem 1.6]{ST} for a more detailed proof).

\bl\label{lemma:method minimal resol}\cite{GS}
Let $A=K\Gamma/I$ be a finite-dimensional $K$-algebra defined by quiver $\Gamma$ and relations $I$. For $i$ in the set of vertices $\Gamma_0$ of $\Gamma$, denote by $e_i$ the corresponding idempotent and by $S_i$ the corresponding simple right $A$-module.  Let  $(P^{\bullet},d^{\bullet})$ be a minimal projective right $A$-module resolution of $A/\rad A$.

Let $(Q^{\bullet},\partial^{\bullet})$ be  a complex  of $A$-$A$-bimodules with $Q^{-1}=A$ and $Q^n=\bigoplus_{i\in Q_0}(Ae_i\ot e_jA)^{\dim\Ext^n_A(S_i,S_j)}$ for $n\pgq 0.$ Assume that $((A/\rad A)\otimes_AQ^{\bullet},\id\ot_A\partial^{\bullet})=(P^{\bullet},d^{\bullet})$. 

Then  $(Q^{\bullet},\partial^{\bullet})$ is a minimal projective resolution of $A$ as an $A$-$A$-bimodule.
\el

\br Note that once a space $Q^2$ satisfying the conditions in the Lemma is found, a minimal set of relations for the ideal $I$ is then given by $(A/\rad{A})^e\ot_{A^e}Q^2$. 
\er

\subsection{Some Lie algebra invariants}

Let $\mg$ be a finite dimensional Lie algebra over $K$ with bracket $[,]$. We briefly recall here a few objects associated to $\mg$ that we will use throughout the paper.

The \textbf{lower central series} of $\mg$, whose $i$-th term is denoted by $\L^i(\mg)$, is defined inductively by 
\[ \L^0(\mg)=\mg,\ \L^1(\mg)=[\mg,\mg]\text{ and }\L^i(\mg)=[\mg,\L^{i-1}(\mg)]\text{ for }i\pgq 2. \] If $\L^i(\mg)=0$ for $i$ large enough, the Lie algebra $\mg$ is \textbf{nilpotent}. 

The \textbf{derived series} of $\mg$,  whose $i$-th term is denoted by $\D^i(\mg)$, is defined inductively by 
\[ \D^0(\mg)=\mg,\ \D^1(\mg)=[\mg,\mg]\text{ and }\D^i(\mg)=[\D^{i-1}(\mg),\D^{i-1}(\mg)]\text{ for }i\pgq 2. \] 

%The \textbf{centre} $\Z(\mg)$ of $\mg$ is the set $\set{x\in \mg\,;\, \forall y\in\mg,\ [x,y]=0}$.

The \textbf{nilradical} of $\mg$ is the maximal nilpotent ideal in $\mg$.

The lower central series, derived series and nilradical are clearly preserved under isomorphisms of Lie algebras.

We now recall the \textbf{Killing form} of $\mg$. This is the bilinear form $\kappa\colon \mg\times\mg\rightarrow K$ defined by 
\[ \kappa(x,y)=\trc([x,[y,-]]). \] If $\mg$ and $\mg'$ are isomorphic Lie algebras, then their Killing forms are equivalent. In particular, they have the same rank.

Finally, we introduce generalised derivations of $\mg$, that were defined in \cite{NH} and that we will use in the proof of Proposition \ref{prop:sd21 constants non zero}. Let $\lambda,$ $\mu$, $\nu$ be three elements in $K$ that are not all zero. A $(\lambda,\mu,\nu)$-\textbf{derivation} of $\mg$ is a linear map $D\colon \mg\rightarrow \mg$ that satisfies 
\[ \lambda D([x,y])=\mu[D(x),y]+\nu[x,D(y)]\text{ for all } x,y,z\text{ in }\mg. \] Let $\der_\mg(\lambda,\mu,\nu)$ denote the space of  $(\lambda,\mu,\nu)$-derivations of $\mg$.

As was mentioned by Novotný and Hrivnák in \cite[Proposition 3.1]{NH},  if $\mg$ and $\mg'$ are isomorphic Lie algebras, then $\der_\mg(\lambda,\mu,\nu)$ and $\der_{\mg'}(\lambda,\mu,\nu)$ are isomorphic vector spaces.

\section{Algebras of dihedral type}\label{sec:dihedral}

The only remaining question in the classification of the algebras of dihedral type up to stable equivalence of Morita type are whether the local algebras $D(1\A)_2^k(d)$ with $d\in\set{0,1}$ in characteristic $2$ are equivalent or not.

Fix an integer $k\pgq 2.$ Consider the local tame symmetric  algebras of dihedral type $\Lambda:=D(1\A)_2^k(d)=K\rep{x,y}/I_d^k$ where
$I_d^k$ is the ideal generated by $\set{x^2-(xy)^k,y^2-d(xy)^k;(xy)^k-(yx)^k;(xy)^kx;(yx)^ky}$ for
$d\in\set{0,1}.$ As we explained in Subsection \ref{subsec:questions:dihedral},  we must determine whether these two algebras are
 equivalent or not. We shall see that the first cohomology group $\HH^1(
D(1\A)_2^k(d))$ enables us to do this.

\bl\label{lemma:proj resol dihedral 1 simple} Consider the sequence of $\Lambda$-$\Lambda$-bimodules \[
Q^2=\bigoplus_{i=0}^2(\Lambda\s{i}\Lambda)\stackrel{\partial^2}{\longrightarrow}Q^1=(\Lambda\s{x}\Lambda)\oplus(\Lambda\s{y}\Lambda)
\stackrel{\partial^1}{\longrightarrow}Q^0=\Lambda\ot \Lambda
\stackrel{\partial^0}{\longrightarrow}\Lambda\rightarrow0\] with the maps determined by 
\begin{align*}
\partial ^1(1\s{\delta}1)&=\delta\ot 1+1\ot \delta\text{ for $\delta\in\set{x,y}$}\\
\partial^2(1\s{0}1)&=x\s{x} 1+1\s{x} x+\sum_{t=0}^{k-1}\left((xy)^t\s{x} y(xy)^{k-1-t}+(xy)^tx\s{y} (xy)^{k-1-t}\right)\\
\partial^2(1\s1 1)&=\sum_{t=0}^{k-1}\left((xy)^t\s{x} y(xy)^{k-1-t}+(xy)^tx\s{y} (xy)^{k-1-t}\right.\\&\qquad\left.+(yx)^t\s{y} x(yx)^{k-1-t}+(yx)^ty\s{x} (yx)^{k-1-t}\right)\\
\partial^2(1\s2 1)&=y\s{y} 1+1\s{y} y+ d\sum_{t=0}^{k-1}\left((yx)^t\s{y} x(yx)^{k-1-t}+(yx)^ty\s{x} (yx)^{k-1-t}\right),
\end{align*}  where the subscripts under the tensor product symbols $\ot$ denote the component of the free
$\Lambda$-$\Lambda$-bimodule $Q^n$ for $n=1,2$.

Then this is the beginning of a minimal projective  $\Lambda$-$\Lambda$-bimodule resolution of $\Lambda.$
\el

\bpf It is easy to check that it is  a complex, and that applying $(\Lambda/\rad\Lambda)\ot_\Lambda?$ gives the beginning of  a minimal projective right $\Lambda$-module resolution of $K=\Lambda/\rad\Lambda$. From this resolution, we may determine $\dim \Ext^n_\Lambda(K,K)$ for $n=0,$ $1$ and $2$ and check that $Q^n$ is the module in Happel's theorem \cite{H}. We then apply Lemma \ref{lemma:method minimal resol}.
\epf

We shall now determine $\HH^1(\Lambda)$. Recall that the centre of $\Lambda$  has dimension $k+3$ and that it is isomorphic to
$\HH^0(\Lambda)=\ker(?\circ \partial^1)$. Therefore the dimension of the image of the map
$(?\circ \partial^1\colon \Hom_{\Lambda-\Lambda}(Q^0,\Lambda)\rightarrow\Hom_{\Lambda-\Lambda}(Q^1,\Lambda))$
is equal to \[\dim
\Hom_{\Lambda-\Lambda}(\Lambda\ot\Lambda,\Lambda)-\dim \HH^0(\Lambda)=4k-(k+3)=3k-3.\]

In order to determine  the kernel of the map
$(?\ot\partial^2\colon \Hom_{\Lambda-\Lambda}(Q^1,\Lambda)\rightarrow\Hom_{\Lambda-\Lambda}(Q^2,\Lambda))$, we note that 
 an element in $\Hom_{\Lambda-\Lambda}(Q^1,\Lambda)=\displaystyle\bigoplus_{i\in\set{x,y}}\Hom_{\Lambda-\Lambda}(\Lambda\s{i}\Lambda,\Lambda)$ is determined by 
\[
f(1\s{i}1)=\sum_{t=0}^{k-1}\left(\lambda_t^{(i)}(xy)^t+\mu_t^{(i)}(yx)^{t+1}+\tau_{t}^{(i)}y(xy)^{t}+\sigma_t^{(i)}x(yx)^t\right) \]
for $i\in\set{x,y}$, where $\lambda_t^{(i)}$, $\mu_t^{(i)},$ $\sigma_t^{(i)}$ and $\tau_t^{(i)}$ are
scalars. Note that $f\circ \partial^2(1\s1 1)=0$ for any $f$.

We then determine the conditions on the coefficients for $f\circ \partial ^2$ to vanish, using standard linear algebra.
We obtain $\dim\ker(?\circ\partial^2)=
\begin{cases}
4k+3&\text{ if $k$ is even and $d=0$}\\
4k+2&\text{ if $k$ is odd and $d=0$}\\&\text{ or if $k$ is even and $d=1$}\\
4k+1&\text{ if $k$ is odd and $d=1$}.
\end{cases}
$ 

Hence we have the following result.

\bp The first cohomology group $\HH^1(D(1\mathcal{A})_2^k(d))$ has dimension $\begin{cases}
k+6-d&\text{ if $k$ is even}\\
k+5-d&\text{ if $k$ is odd.}
\end{cases}$
\ep 

\bc\label{cor:local dihedral} There is no stable equivalence  of
Morita type between the algebras  $D(1\mathcal{A})_2^k(0)$ and $D(1\mathcal{A})_2^k(1)$.
\ec

\br This completes the classification of the algebras of dihedral type up to stable equivalence of Morita type.
\er

\section{Algebras of semi-dihedral type}\label{sec:semidihedral}

As we mentioned in Subsection \ref{subsec:questions:semidihedral}, the classification is complete for algebras of semi-dihedral type with three simple modules. We shall start with the local algebras.

\subsection{Local algebras of semi-dihedral type}

\subsubsection{Dimension of the first Hochschild cohomology group}

\sloppy We assume here that the field $K$ has characteristic $2.$ Fix an integer $k\pgq2.$ For $(c,d)\in K^2$, let $I^k(c,d)$ be the ideal in $K\rep{x,y}$
generated by $\set{(xy)^k+(yx)^k;x^2+(yx)^{k-1}y+c(yx)^k;y^2+d(xy)^k;(xy)^kx)}.$ For any local tame
symmetric algebra of semi-dihedral type $\Lambda$, there  is a stable
equivalence of Morita type with one of the algebras $SD(1\A)_1^k=K\rep{x,y}/I^k(0,0)$ and
$SD(1\A)_2^k(c,d)=K\rep{x,y}/I^k(c,d)$ for $(c,d)\in K^2\backslash\set{(0,0)}$ (see
\cite[Theorem 7.1]{ZZ}). However, for a fixed $k,$ it is not known whether these algebras are stably equivalent of Morita type or not.

\sloppy Using isomorphisms of the form $(x,y)\mapsto (\lambda x,\mu y)$, we can assume that $(c,d)\in\set{(1,0);(c,1);\text{ with } c\in K}$.  Note that in all these algebras, we have the following identities: 
\[ xy^2=0=y^2x;\ y(xy)^k=0;\ x^2y=0=yx^2;\ x^3=(xy)^k=(yx)^k;\ x^4=0;\  y^3=0. \]

The aim of this section is to prove the following theorem.

\bt\label{thm:classification sd local}
The algebras $SD(1\A)_1^k$, $SD(1\A)_2^k(1,0)$, $SD(1\A)_2^k(0,1)$ and
$SD(1\A)_2^k(c,1)$ for $c\neq 0$ are in  four different stable equivalence of Morita type
classes.
\et

\br    We cannot say whether $SD(1\A)_2^k(c,1)$ and $SD(1\A)_2^k(c',1)$ for $c\neq c'$
non-zero are stably equivalent of Morita type or not. 
\er

In the rest of the section, $\Lambda$ is one of the algebras $SD(1\A)_1^k$ or $SD(1\A)_2^k(c,d)$.

\bl\label{lemma:proj resol semi-dihedral 1 simple} Consider the sequence of $\Lambda$-$\Lambda$-bimodules \[
Q^2=\bigoplus_{i=0}^1(\Lambda\s{i}\Lambda)\stackrel{\partial^2}{\longrightarrow}Q^1=(\Lambda\s{x}\Lambda)\oplus(\Lambda\s{y}\Lambda)
\stackrel{\partial^1}{\longrightarrow}Q^0=\Lambda\ot \Lambda
\stackrel{\partial^0}{\longrightarrow}\Lambda\rightarrow0\] with the maps determined by 
\begin{align*}
\partial ^1(1&\s \delta 1)=\delta\ot 1+1\ot \delta\text{ for }\delta\in\set{x,y}\\
\partial^2(1&\s{0}1)=x\s{x}1+1\s{x}x+\sum_{t=0}^{k-2}(yx)^ty\s{x}y(xy)^{k-2-t}\\&+\sum_{t=0}^{k-1}\left((yx)^t\s{y}(xy)^{k-1-t}+c(yx)^ty\s{x}(yx)^{k-1-t}+c(yx)^{t}\s{y}x(yx)^{k-1-t}\right)\\
\partial^2(1&\s1 1)=y\s{y}1+1\s{y}y+d\sum_{t=0}^{k-1}\left((xy)^tx\s{y}(xy)^{k-1-t}+(xy)^t\s{x}y(xy)^{k-1-t}\right),
\end{align*} where the subscripts on the tensor product symbols $\ot$ denote the component of the free
$\Lambda$-$\Lambda$-bimodule $Q^n$ for $n=1,2$.

Then this is the beginning of a minimal projective  $\Lambda$-$\Lambda$-bimodule resolution of $\Lambda.$
\el

\bpf The proof is the same as that of  Lemma \ref{lemma:proj resol dihedral 1 simple}.
\epf

Using this resolution,  we may compute the Hochschild cohomology groups. 

As in the case of the dihedral algebras $D(1\A)_2^k(d)$, we have $\dim \im(?\circ \partial ^1)=3k-3$ (we give a generating set explicitly in the proof of Lemma \ref{lemma:basis sd1}). Moreover, it is easy to check that 
 $\dim\ker(?\circ \partial ^2)=
\begin{cases}
4k+3&\text{ if $k$ is even and $d=0$ or if $k$ is odd and $c=0=d$,}\\
4k+2&\text{ if $k$ is even and $d\neq0$ or if $k$ is odd, $c\neq 0$ and $d=0$,}\\
4k+1&\text{ if $k$ is odd and $d\neq 0$.}\\
\end{cases}
$  

Therefore we get the following dimensions for the first Hochschild cohomology group.

\bp\label{prop:HH1 semi-dihedral 1 simple} Let $\Lambda$ be one of the algebras $SD(1\A)_1^k$ or $SD(1\A)_2^k(c,d)$. Then
\[ \dim\HH^1(\Lambda)=\begin{cases}
k+6&\text{ if $k$ is even and $d=0$ or if $k$ is odd and $c=0=d$,}\\
k+5&\text{ if $k$ is even and $d\neq0$ or if $k$ is odd, $c\neq 0$ and $d=0$,}\\
k+4&\text{ if $k$ is odd and $d\neq 0$.}\\
\end{cases}
 \]
\ep

\bc\label{cor:HH1 semi-dihedral 1 simple} For any $k\pgq 2$,  an algebra in the set $\set{SD(1\A)_1^k;SD(1\A)_2^k(1,0)}$ is  not  stably equivalent of Morita type to an algebra in the set $\set{SD(1\A)_2^k(0,1);SD(1\A)_2^k(c,1), c\neq 0}$.

Moreover, if $k$ is odd there is no stable equivalence of Morita type between the algebras $SD(1\A)_1^k$ and $SD(1\A)_2^k(1,0)$.
\ec

\subsubsection{Lie algebra structure on $\HH^1(\Lambda)$}\label{subsection:Lie semidihedral}

We shall now improve on Corollary \ref{cor:HH1 semi-dihedral 1 simple} using the Lie algebra structure on $\HH^1(\Lambda)$ where $\Lambda$ is one of the algebras $SD_1^k(1\A)$ or
$SD_2^k(1\A)(c,d).$ Let $\Gamma$ be a quiver of type $1\A$,  with arrows $x$ and $y$. 
Then the local tame symmetric algebras of semi-dihedral type
 may be defined as $K\Gamma/I^k(c,d)$ for $(c,d)\in K^2$ or, as we mentioned above, $(c,d)\in\set{(0,0),(1,0),(c,1);c\in K}.$ 

\br It is possible (though laborious) in this case to compute $\dim\HH^n(\Lambda)$ for all $n\pgq 0$ (the case $c=0=d$ may be found in \cite{G10}). However, these dimensions do not give any more information than  $\dim\HH^1(\Lambda)$.
\er

We  have $\Hom_{\Lambda-\Lambda}((\Lambda\s x\Lambda)\oplus(\Lambda\s y\Lambda),\Lambda)\cong\Hom_{\Lambda-\Lambda}(\Lambda\ot_KK\Gamma_1\ot_K\Lambda,\Lambda)\cong\Hom_K(K\Gamma_1,\Lambda)$ via the correspondence $f\leftrightarrow g$ given by $f(1\s x1)=g(x)$ and $f(1\s y1)=g(y).$ We shall often identify $g\in \Hom_K(K\Gamma_1,\Lambda)$ with the pair $(g(x),g(y))$.

Define the following elements in $\Hom_K(K\Gamma_1,\Lambda):$
\begin{align*}
&\varphi_t=(x(yx)^t,0) \text{ for } 0\ppq t\ppq k-1, && \theta_0=(1+cx,cy+d(yx)^{k-1}),\\&\theta_1=(y(xy)^{k-1},0),&&
\theta_{-1}=(0,x(yx)^{k-1}),\\&\theta_2=((xy)^k,0),&&\theta_{-2}=(0,(xy)^k),\\
&\omega=(y(xy)^{k-2}+c(yx)^{k-1},1),&&\chi=(0,y).
\end{align*}

Set $\B=\set{\varphi_t, 1\ppq t\ppq k-1;\theta_1;\theta_{-1};\theta_2;\theta_{-2};\theta_0}\subset\Hom_K(K\Gamma_1,\Lambda)$ and \[\B'=
\begin{cases}
\set{\omega;\varphi_0}&\text{ if $k$ is odd and $c=0=d$},\\
\set{\omega;\chi}&\text{ if $k$ is even and $d=0$},\\
\set{\chi}&\text{ if $k$ is even and $d\neq 0$},\\
\set{\omega}&\text{ if $k$ is odd, $c\neq 0$ and $d=0$},\\
\varnothing&\text{ if $k$ is odd and $d\neq 0$}.
\end{cases}\] 

\bl \label{lemma:basis sd1}
With the notation above, $\B\cup\B'$ is a set of cocycle representatives of a basis of $\HH^1(\Lambda)$.
\el

\bpf The fact that the elements in $\B\cup\B'$ are cocycles can be checked easily (recall that $\car(K)=2$).

Moreover,  the classes of the cocycles \[ \C:=\set{(x(yx)^t,y(xy)^t);((xy)^t+(yx)^t,0);(0,(xy)^t+(yx)^t);1\ppq t\ppq k-1} \] form a basis of $\im(?\circ \partial ^1)$ (this basis will be useful when computing Lie brackets).

It is then straightforward to check that the cochains in $\B\cup\B'\cup\C$ are linearly independent, and the result follows, using the dimension of $\HH^1(\Lambda)$ obtained previously.
\epf

As described in Subsection \ref{sec:general description Lie algebra}, we  transport the usual Lie bracket on $\HH^1(\Lambda)$
defined using the Bar resolution to a Lie bracket on $\HH^1(\Lambda)$ defined using the minimal projective
resolution. Note that we can identify $Q^2$ with $\Lambda\ot_{K\Gamma_0}KZ\ot_{K\Gamma_0}\Lambda$ with $Z=\set{x^2-(yx)^{k-1}y+c(yx)^k;y^2-d(xy)^k}.$

\bl\label{lemma:sd local brackets} We use the same notation for a cocycle and for its cohomology class. The (potentially) non-zero brackets of basis elements are the following:
\begin{alignat*}{3}
&[\varphi_t,\varphi_{t'}]=(t+t')\varphi_{t+t'};&\quad& [\theta_{-2};\chi]=\theta_{-2};&\quad& [\theta_{-1},\varphi_0]=\theta_{-1};\\
&[\varphi_t,\theta_0]=
\begin{cases}
\theta_0\text{ if $t=0$}\\d(k-1)\theta_{-2}\text{ if $t=1$}\\0\text{ if $t>1$};
\end{cases}&\quad&  [\theta_{-2},\varphi_0]=\theta_{-2};&\quad& [\theta_1;\varphi_0]=\theta_1;\\
& [\varphi_t,\omega]=
\begin{cases}
(k-1)(\theta_1+c\theta_2)\text{ if $t=1$}\\0\text{ if $t\neq 1$};
\end{cases}&\quad& [\theta_{-2},\theta_0]=\varphi_{k-1}+c\theta_{-2};&\quad& [\theta_0,\omega]=c\omega;\\
&[\varphi_t,\chi]=t\varphi_t;&\quad&[\theta_2,\omega]=\varphi_{k-1} ;&\quad& [\omega,\chi]=\omega.\\
& [\theta_{-2},\omega]=
\begin{cases}
\theta_{-1}\text{ if $k>2$}\\\theta_{-1}+\theta_2\text{ if $k=2$};
\end{cases}&\quad& [\theta_2,\theta_0]=\theta_1+c\theta_2;
\end{alignat*}
\el 

\bpf We refer to Lemma \ref{lemma:sd 2 simples brackets} for an example (in a non-local case) of the computation of a Lie bracket, the method here is similar.
\epf

We then compute the first two terms in the lower central series. Since they give no
new information when $d=0,$ we only give the results for $d\neq0.$

\bp We keep the notation above and assume that $d\neq 0.$ Then a basis of $\L^1(\HH^1(SD(1\A)_2^k(c,d)))$
is given by:
\begin{enumerate}[(a)]
\item $\set{\varphi_{2p+1},1\ppq p\ppq \frac{k-3}{2};\theta_1+c\theta_2;\varphi_{k-1}+c\theta_{-2}}$
  if $k$ is odd (there are no $\varphi_t$ if $k=3$);
\item $\set{\varphi_{2p+1},0\ppq p\ppq \frac{k-2}{2};\theta_1+c\theta_2;\theta_{-2}}$ if $k$ is even.
\end{enumerate}

Moreover, $\L^2(\HH^1(SD(1\A)_2^k(c,d)))$ is generated by the following set:
\begin{enumerate}[(a)]
\item $\set{\varphi_{2p+1},2\ppq p\ppq
    \frac{k-3}{2};c(\theta_1+c\theta_2);c(\varphi_{k-1}+c\theta_{-2})}$ if $k$ is odd (there are no
  $\varphi_t$ if $k=3$ or $k=5$);
\item $\set{\varphi_{2p+1},0\ppq p\ppq \frac{k-2}{2};c(\theta_1+c\theta_2);\theta_{-2}}$ if $k$ is even.
\end{enumerate} where $\delta$ is the Kronecker symbol.
In particular, $\dim \L^2(\HH^1(SD(1\A)_2^k(c,d)))=
\begin{cases}
\frac{k-1}{2}-2\delta_{c,0}+\delta_{k,3}&\text{ if $k$ is odd},\\
\frac{k}{2}+2-\delta_{c,0}&\text{ if $k$ is even.}
\end{cases}
$
\ep

\bc 
There is no stable equivalence of Morita type between the algebras   $SD(1\A)_2^k(0,1)$ and $SD(1\A)_2^k(c,1)$ for $c\neq 0$.
\ec

In order to complete the proof of Theorem \ref{thm:classification sd local}, we must prove the following result.

\bp\label{prop:sd 2 simples d=0}
Assume that $k$ is even. Then there is no stable equivalence of Morita type between the algebras   $SD(1\A)_1^k$ and  $SD(1\A)_2^k(1,0)$. 
\ep

\bpf Let $\mg$ be the Lie algebra $\HH^1(SD(1\A)_1^k)$ and $\mg'$ be the Lie algebra $\HH^1(SD(1\A)_2^k(1,0))$. These Lie algebras are not nilpotent (indeed, since $[\varphi_1,\chi]=\varphi_1$, it follows that $\varphi_1$ is in all the terms of the lower central series for both Lie algebras). 

Consider the subspace $I$ of $\mg$ generated by $\set{\varphi_t,1\ppq t\ppq k-1;\theta_0;\theta_{-1};\theta_1;\theta_{-2};\theta_2;\omega}$ and the subset $I'$ of $\mg'$ generated by $\set{\varphi_t,1\ppq t\ppq k-1;\theta_{-1};\theta_1;\theta_{-2};\theta_2;\omega}$. They are Lie ideals.

Moreover, the lower central series of $I$ is given by $\L^1(I)=\spn{\varphi_t,3\ppq t\ppq k-1;\theta_0;\theta_{-1};\theta_1}$ and $\L^i(I)=\spn{\varphi_t,{2^{i+1}-1}\ppq t\ppq k-1}$ if $i\pgq 2$, so that it vanishes eventually and $I$ is nilpotent. Similarly,  the lower central series of $I'$ is given by $\L^1(I')=\spn{\varphi_t,3\ppq t\ppq k-1;\theta_{-1};\theta_1+\theta_2}$ and $\L^i(I')=\spn{\varphi_t,{2^{i+1}-1}\ppq t\ppq k-1}$ if $i\pgq 2$, so that it vanishes eventually and $I'$ is nilpotent.

Since $\dim I=\dim \mg-1$ and $\mg$ is not nilpotent, $I$ is the nilradical of $\mg$.

We now prove that $I'$ is the nilradical of $\mg'$. Assume for a contradiction that it is not. Then
it follows that there is a non-zero element in $\mg'$, that we can choose of the form
$u=\lambda\chi+\mu\theta_0$, such that the subspace $J$ generated by $I'$ and $u$ is a nilpotent
ideal. Since $[u,\omega]=(\lambda+\mu )\omega$ and $J$ is nilpotent, we must have $\lambda+\mu =0$
(otherwise $\omega$ would be in all the $\L^i(J)$). Therefore we may assume that $u=\chi+\theta_0$. We have $[u,\theta_1+\theta_2]=\theta_1+\theta_2$ so that $\theta_1+\theta_2\in \L^i(J)$ for all $i$, a contradiction. 
Therefore $I'$ is the nilradical of $\mg'$.

It follows that the nilradicals of  $\HH^1(SD(1\A)_1^k)$ and $\HH^1(SD(1\A)_2^k(1,0))$ have different dimensions, and hence that  $\HH^1(SD(1\A)_1^k)$ and $\HH^1(SD(1\A)_2^k(1,0))$ are not isomorphic Lie algebras.
\epf

\br It can be noted that when $k$ is odd and $cc'\neq 0$, the Lie algebras $\HH^1(SD(1\A)_2^k(c,1))$ and  $\HH^1(SD(1\A)_2^k(c',1))$ are isomorphic. Indeed, if $\set{\varphi_t, 1\ppq t\ppq k-1;\theta_1;\theta_{-1};\theta_2;\theta_{-2};\theta_0}$ is a basis of  $\HH^1(SD(1\A)_2^k(c,1))$ and  $\set{\varphi_t', 1\ppq t\ppq k-1;\theta_1';\theta_{-1}';\theta_2';\theta_{-2}';\theta_0'}$ is a basis of  $\HH^1(SD(1\A)_2^k(c',1))$, the isomorphism is defined by 
\begin{align*}
&\varphi_t\mapsto\varphi_t',&&\theta_1\mapsto \theta_1',&&\theta_{-1}\mapsto \theta_{-1}',\\
&\theta_0\mapsto \frac{c}{c'}\theta_0',&&\theta_{2}\mapsto \frac{c'}{c}\theta_{2}',&&\theta_{-2}\mapsto \frac{c'}{c}\theta_{-2}'.
\end{align*} In the remaining unresolved cases, we do not know whether the first Hochschild cohomology groups are isomorphic or not.
\er

\subsection{Algebras of semi-dihedral type with two simple modules}\label{subsec:sd2}

We have defined the algebras  $SD(2\B)_1^{k,s}(c)$ and $SD(2\B)_2^{k,s}(c)$ of semi-dihedral type with two simple modules in Section \ref{sec:questions}. 
Note that when  $k\pgq 2$, the ideal of relations for  $SD(2\B)_2^{k,2}(c)$ is not admissible; a definition with an admissible ideal can be obtained by removing the loop $\eta$ and adapting the relations.

We shall use the Lie algebra structure of the first Hochschild cohomology group to improve on the results in \cite{ZZ}, and to give a partial answer to the question of  whether the algebras $SD(2\B)_1^{k,s}(c)$ and $SD(2\B)_2^{k',s'}(c')$ are stably  equivalent of Morita type or not. 

The main result of this section is the following.

\bt\label{thm:sd} Let $k,k',s,s'$ be integers with $\set{k,s}=\set{k',s'}$ and let $c$ be an element in $\set{0,1}$.
\begin{enumerate}[(1)]
\item Assume that $\car(K)=2.$ Then, for $a\in\set{1,2}$, the algebras $SD(2\B)_a^{k,s}(0)$ and $SD(2\B)_a^{k',s'}(1)$ are not stably  equivalent of Morita type.
\item In each of the following cases,  the algebras $SD(2\B)_1^{k,s}(c)$ and $SD(2\B)_2^{k',s'}(c')$  are not stably   equivalent of Morita type:
\begin{enumerate}[(i)]
\item $\car(K)=2$ and $ks$ is even;
\item  $\car(K)=2$, $ks$ is odd and $(c,c')\neq (0,0)$;
\item $\car(K)=3$; 
\item $\car(K)\neq 2,3$, $ks=0$ in $K$;
\item $\car(K)\neq 2,3$, $ks\neq 0$ in $K$ and $k+s-2ks=0$ in $K$.
\item $\car(K)\neq 2,3$, $\lambda:=3^{-1}2ks\neq 0$ in $K$, $\mu:=2ks-k-s\neq 0$ in $K$, and the following subsets of $K$ are not equal:
\begin{align*}
&\set{s\lambda^{-1},2s\lambda^{-1},k\lambda^{-1},2k\lambda^{-1},(s\lambda^{-1})^{-1},(2s\lambda^{-1})^{-1},(k\lambda^{-1})^{-1},(2k\lambda^{-1})^{-1}}\\
and &\set{s\mu^{-1},2s\mu^{-1},k\mu^{-1},2k\mu^{-1},(s\mu^{-1})^{-1},(2s\mu^{-1})^{-1},(k\mu^{-1})^{-1},(2k\mu^{-1})^{-1}}.
\end{align*}
\end{enumerate}
\end{enumerate}
\et

The remainder of Section \ref{subsec:sd2} is devoted to the proof of this result.

\subsubsection{The first Hochschild cohomology group of $SD(2\B)_1^{k,s}(c)$}
\allowdisplaybreaks
Let $\Lambda$ be the algebra $SD(2\B)_1^{k,s}(c)$ and let $\Gamma$ be the quiver of type $2\B$. Let $e_1$ and $e_2$ denote the idempotents in $\Lambda$ corresponding to the vertices.

\bl\label{lemma:proj resol semidihedral 2 simples} Define a sequence of $\Lambda$-$\Lambda$-bimodules $
Q^2 \stackrel{\partial^2}{\longrightarrow}Q^1
\stackrel{\partial^1}{\longrightarrow}Q^0
\stackrel{\partial^0}{\longrightarrow}\Lambda\rightarrow0$ as follows. The modules $Q^n$ are given by 
\begin{align*}
&Q^2=(\Lambda e_1\ot e_1\Lambda )\oplus (\Lambda e_1\ot e_2\Lambda )\oplus (\Lambda e_2\ot e_1\Lambda )\oplus (\Lambda e_2\s1 e_2\Lambda )\oplus (\Lambda e_2\s2 e_2\Lambda ),\\ &Q^1=\bigoplus_{i,j=1}^2\Lambda e_i\ot e_j\Lambda \\
&Q^0=\bigoplus_{i=1}^2\Lambda e_i\ot e_i\Lambda ,
\end{align*} where the subscripts on the tensor product symbols $\ot$ denote the component of the free
$\Lambda $-$\Lambda $-bimodule $Q^2$. The map  $\partial ^0$ is multiplication and the other maps are determined by \begin{align*}
&\partial ^1(e_{i(\delta)}\ot e_{t(\delta)})=\delta\ot e_{t(\delta)}-e_{i(\delta)}\ot \delta\quad\text{ for }\delta\in\Gamma_1\text{ with origin $i(\delta)$ and endpoint $t(\delta)$}\\
&\partial^2(e_1\ot e_1)=e_1\ot \alpha+\alpha\ot e_1-\sum_{t=0}^{k-1}(\beta\gamma\alpha)^t\left(\beta\ot e_1+e_1\ot \gamma\right)(\alpha\beta\gamma)^{k-1-t}\\&\qquad\qquad-\sum_{t=0}^{k-2}(\beta\gamma\alpha)^t\beta\gamma\ot \beta\gamma(\alpha\beta\gamma)^{k-2-t}\\&\qquad\qquad+c\sum_{t=0}^{k-1}(\beta\gamma\alpha)^t\left(e_1\ot \gamma\alpha+\beta\ot \alpha+\beta\gamma\ot e_1\right)(\beta\gamma\alpha)^{k-1-t}\\
&\partial^2(e_1\ot e_2)=e_1\ot\eta+\beta\ot e_2\\
&\partial^2(e_2\ot e_1)=e_2\ot\gamma+\eta\ot e_1\\
&\partial^2(e_2\s1 e_2)=e_2\ot\beta+\gamma\ot e_1\\
&\partial^2(e_2\s2 e_2)=\sum_{ r=0}^{s-1}\eta^ r\ot\eta^{s-1- r}-\sum_{t=0}^{k-1}(\gamma\alpha\beta)^t\left(e_2\ot\alpha\beta+\gamma\ot\beta+\gamma\alpha\ot e_2\right)(\gamma\alpha\beta)^{k-1-t}.
\end{align*}

Then this sequence is the beginning of a minimal projective  $\Lambda $-$\Lambda $-bimodule resolution of $\Lambda .$ 
\el

\bpf  It is easy to check that it is  a complex, and that applying $S_i\ot_\Lambda ?$ gives the beginning of  a minimal projective right $\Lambda $-module resolution of the simple module $S_i$ for $i=1,2$. From these resolutions, we may determine $\dim \Ext^n_\Lambda (S_i,S_j)$ for $n=0,$ $1$ and $2$ and $i,j=1,2$ and check that $Q^n$ is the module in Happel's theorem \cite{H}. Noting that $\Lambda /\rad \Lambda =S_1\oplus S_2$ as a right $\Lambda $-module, we then apply Lemma \ref{lemma:method minimal resol}.
\epf

\br We can identify $Q^1$ with $\Lambda\ot_{K\Gamma_0} K\Gamma_1\ot _{K\Gamma_0} \Lambda$ via $a\ot \delta\ot a'\mapsto ae_{i\delta}\ot e_{t(\delta)}a'$ and similarly $Q^2$ with $\Lambda\ot _{K\Gamma_0}KZ\ot_{K\Gamma_0}\Lambda$ where $Z=\set{\alpha^2-\beta\gamma(\alpha\beta\gamma)^{k-1}-c(\alpha\beta\gamma)^k;\beta\eta;\eta\gamma;\gamma\beta;\eta^s-(\gamma\alpha\beta)^k}.$
\er

Using the resolution above, we may now compute the dimension of $\HH^1(\Lambda )$. The proof is straightforward and is omitted.

\bp Let  $\Lambda $ be the algebra $SD(2\B)_1^{k,s}(c)$. Then 
\[ \dim\HH^1(\Lambda )=
\begin{cases}
k+s+3&\text{ if }\car(K)=2 \text{ and $k$ and $s$ are both even } \\
k+s+2&\text{ if }\car(K)=2 \text{ and $k$ and $s$ are not both even}\\&\qquad\qquad\text{ and $ksc=0$ in $K$} \\
k+s+1&\text{ if }\car(K)=2 \text{ and $ksc\neq 0$ in $K$} \\
k+s+2&\text{ if }\car(K)=3 \text{ and $k$ and $s$ are both $0$ in $K$} \\
k+s+1&\text{ if }\car(K)=3 \text{ and $k$ and $s$ are not both $0$ in $K$}\\&\qquad\qquad\text{  but $ks=0$ in $K$} \\
k+s&\text{ if }\car(K)=3 \text{ and $ks\neq 0$ in $K$} \\
k+s+1&\text{ if }\car(K)\neq 2,3 \text{ and $k$ and $s$ are both $0$ in $K$} \\
k+s&\text{ if }\car(K)\neq 2,3 \text{ and $k$ and $s$ are not both $0$ in $K$.} 
\end{cases}
 \]
\ep

\bc\label{cor:sd1 car2} If $\car(K)=2$ and $k$ and $s$ are both odd and if $\set{k',s'}=\set{k,s}$, then there is no stable equivalence of Morita type between $SD(2\B)^{k,s}_1(0)$ and $SD(2\B)^{k',s'}_1(1).$
\ec

In order to go further, we now consider the Lie algebra structure of $\HH^1(\Lambda )$.  In the sequel, we identify a morphism $f\in\Hom_{\Lambda -\Lambda }(Q^1,\Lambda )$ with $g\in\Hom_{K\Gamma_0-K\Gamma_0}(K\Gamma_1,\Lambda )$ such that $g(\delta)=f(e_{i(\delta)}\ot e_{t(\delta)})$ for all $\delta\in\Gamma_1$, and with the quadruple $(g(\alpha),g(\beta),g(\gamma),g(\eta)).$

\medskip

\paragraph{\bfseries First assume that $\car(K)=2$} We start with a basis for  $\HH^1(SD(2\B)_1^{k,s}(c))$. 

\bl\label{lemma:sd1 car2 basis} We define cocycles in $\Hom_{K\Gamma_0-K\Gamma_0}(K\Gamma_1,\Lambda )$ as follows.
\begin{align*}
\varphi_t&=(\alpha(\beta\gamma\alpha)^t,0,0,0)\quad\text{ for } 1\ppq t\ppq k-1\\
\theta_ r&=(0,0,0,\eta^{ r+1})\quad\text{ for } 1\ppq r\ppq s-1\\
\psi&=((\alpha\beta\gamma)^k,0,0,0)\\
\chi&=(e_1+c\alpha,c\beta,0,0)\\
\omega&=((\beta\gamma\alpha)^{k-1}\beta\gamma+c(\alpha\beta\gamma)^k,0,0,0)
\\
\varphi_0&=(0,\beta,0,0)\text{ and }\theta_0=(0,0,0,\eta)\text{ if $k$ and $s$ are both even,}\\
\zeta_1&=(0,s\beta,0,k\eta)\text{ if $k+s$ is odd}\\
\zeta_0&=(\alpha,0,0,\eta)\text{ if $k$ and $s$ are both odd and $c=0$.}\\
\end{align*} The cohomology classes of the cocycles in  $\B\cup\B'$ with   $\B=\set{\varphi_t;\theta_ r;\psi;\chi;\omega;1\ppq t\ppq k-1;1\ppq r\ppq s-1}$ and $\B'=
\begin{cases}
\set{\varphi_0;\theta_0}&\text{ if $k$ and $s$ are both even}\\
\set{\zeta_1}&\text{ if $k+s$ is odd}\\
\set{\zeta_0}&\text{ if $k$ and $s$ are both odd and $c=0$}\\
\varnothing&\text{ otherwise}.
\end{cases}
$ \\
form a basis of $\HH^1(SD(2\B)_1^{k,s}(c))$.
\el

\bpf This is proved by computing explicitly $\ker(?\circ \partial ^2)$ and $\im(?\circ \partial ^1)$. We omit the details, but we give the following basis of $\im(?\circ \partial ^1)$, which is useful when computing brackets of elements in $\HH^1(\Lambda)$;
\begin{align*}
&(\alpha(\beta\gamma\alpha)^t,-\beta(\gamma\alpha\beta)^t,0,0),&&(\alpha(\beta\gamma\alpha)^t,0,-\gamma(\alpha\beta\gamma)^t,0),&&((\alpha\beta\gamma)^t-(\beta\gamma\alpha)^t,0,0,0)\\
&(0,\alpha\beta(\gamma\alpha\beta)^t,-\gamma\alpha(\beta\gamma\alpha)^t,0),&&(0,\alpha\beta,-\gamma\alpha,0)&&(0,\beta,-\gamma,0)
\end{align*} with $1\ppq t\ppq k-1.$
\epf

We may now compute the brackets of these basis elements.

\bl\label{lemma:sd 2 simples brackets} We use the notation in the previous lemma. We describe the (potentially) non-zero brackets in $\HH^1(SD(2\B)_1^{k,s}(c))$. 

For all $k $ and $s$ we have
\begin{align*}
&[\varphi_t,\varphi_{t'}]=(t+t')\varphi_{t+t'}\rlap{$ \text{ if }0\ppq t,t'\ppq k-1,\ t+t'\ppq k-1,\ $}\\
&[\theta_r,\theta_{r'}]=(r+r')\theta_{r+r'} \rlap{$  \text{ if }0\ppq r,r'\ppq s-1, \ r+r'\ppq s-1, \ $}\\
&[\psi,\chi]=\omega&
&[\chi,\omega]=c\omega.
\end{align*} 

If moreover $k+s$ is odd, then 
\[ [\varphi_t,\zeta_1]=ts\varphi_t\qquad\text{ and }\qquad [\theta_r,\zeta_1]=rk\theta_r. \]

If instead $k$ and $s$ are both odd and $c=0$, then 
\begin{align*}
 [\varphi_t,\zeta_0]&=t\varphi_t,&[\chi,\zeta_0]&=\chi,& [\theta_r,\zeta_0]&=r\theta_r,&[\omega,\zeta_0]&=\omega.
\end{align*}
\el

\bpf In order to illustrate the method, let us determine the bracket $[\chi,\omega]$. We view $\chi$ and $\omega$ as maps in $\Hom_{K\Gamma_0-K\Gamma_0}(K\Gamma_1,\Lambda ).$

First, for every $\delta\in\Gamma_1$, we replace each instance of $\delta$ in $\omega(\alpha)$ by $\chi(\delta)$, that is, we replace every $\alpha$ in turn with $e_1+c\alpha$ and every $\beta$ by $c\beta$, and we add the results. Since $\car(K)=2$, we get 
\[ (k-1)c(\beta\gamma\alpha)^{k-1}\beta\gamma+kc^2(\alpha\beta\gamma)^k+ce_1\beta\gamma(\alpha\beta\gamma)^k+kc(\beta\gamma\alpha)^{k-1}\beta\gamma +kc^2(\alpha\beta\gamma)^k=0.\] We apply the same procedure to $\omega(\beta)=0,$ $\omega(\gamma)=0$ and $\omega(\eta)=0$, and we obtain $0$ in all cases.

Next, we exchange the roles of $\chi$ and $\omega$. We replace  each instance of $\delta$ in $\chi(\alpha)$ and $\chi(\beta)$ by $\omega(\delta)$. We get 
\begin{align*}
\alpha\mapsto&c(\beta\gamma\alpha)^{k-1}\beta\gamma+c^2(\alpha\beta\gamma)^k\\
\beta\mapsto&0
\end{align*} and of course $\gamma$ and $\eta$ are sent to $0.$ 

Finally, we subtract the two quantities, which gives the map
\begin{align*}
\alpha\mapsto&c(\beta\gamma\alpha)^{k-1}\beta\gamma+c^2(\alpha\beta\gamma)^k\\
\beta\mapsto&0\\
\gamma\mapsto&0\\
\eta\mapsto&0,
\end{align*} that is, $c\omega.$ Therefore, $[\chi,\omega]=c\omega.$

The other brackets are computed in the same way. Note that we work modulo $\im(?\circ \partial ^1).$
\epf

It can be noted that the Lie algebras $\HH^1(SD(2\B)_1^{k,s}(c))$ and  $\HH^1(SD(2\B)_1^{s,k}(c))$ for $c\in \set{0,1}$ (same $c$) are isomorphic.  Indeed, if we consider the basis of  $\HH^1(SD(2\B)_1^{k,s}(c))$ given in Lemma \ref{lemma:sd1 car2 basis} and the similar basis contained in  $\set{\varphi_r',0\ppq r\ppq s-1,\theta_t',0\ppq t\ppq k-1,\psi',\chi',\omega',\zeta_0',\zeta_1'}$ of  $\HH^1(SD(2\B)_1^{s,k}(c))$, the isomorphism is given by 
\[ \varphi_t\mapsto \theta_t',\quad   \theta_r\mapsto \varphi_r',\quad   \psi\mapsto \psi',\quad   \chi\mapsto \chi',\quad   \omega\mapsto \omega',\quad   \zeta_0\mapsto \zeta_0', \quad \zeta_1\mapsto \zeta_1' \] on the elements that are actually present in each case. Therefore the Lie algebra structure of $\HH^1(\Lambda)$ does not help to separate the pairs of parameters $(k,s)$ and $(s,k)$.

We already know from Corollary \ref{cor:sd1 car2} that if $k$ and $s$ are both odd and $\set{k,s}=\set{k',s'}$, then  $SD(2\B)_1^{k,s}(0)$ and $SD(2\B)_1^{k',s'}(1)$ are not stably equivalent of Morita type.

If $k+s$ is odd,  the second term in the lower central series has dimension $\frac{k+s-5}{2}+c+\delta_{k,3}+2\delta_{k,1}$ (if $k$ is odd and $s$ is even, it is spanned by $\set{\varphi_{2p+1};\theta_{2q+1};c\omega;2\ppq p\ppq \frac{k-3}{2},0\ppq q\ppq \frac{s}{2}-1}$), therefore it follows that  $\HH^1(SD(2\B)_1^{k,s}(0))$ and  $\HH^1(SD(2\B)_1^{k',s'}(1))$, with $\set{k,s}=\set{k',s'}$, are not isomorphic Lie algebras, and hence that $SD(2\B)_1^{k,s}(0)$ and $SD(2\B)_1^{k',s'}(1)$ are not stably equivalent of Morita type. 

Similarly, if $k$ and $s$ are both even,  the second term in the lower central series of  $\HH^1(SD(2\B)_1^{k,s}(c))$ is spanned by $\set{\varphi_{2p+1};\theta_{2q+1};c\omega;0\ppq p\ppq \frac{k}{2}-1,0\ppq q\ppq \frac{s}{2}-1}$and  has dimension $\frac{k+s}{2}+c$, therefore  $SD(2\B)_1^{k,s}(0)$ and $SD(2\B)_1^{k',s'}(1)$ are not stably equivalent of Morita type.

We have therefore proved the following result.

\bp\label{prop:sd1 car2} Assume that $\car(K)=2.$ Then the algebras $SD(2\B)_1^{k,s}(0)$ and $SD(2\B)_1^{k',s'}(1)$, with   $\set{k,s}=\set{k',s'}$, are not stably  equivalent of Morita type. 
\ep

\medskip

\paragraph{\bfseries Next assume that $\car(K)=3$} The Lie algebra structure of $\HH^1(SD(2\B)_1^{k,s}(c))$ is determined in the following lemma.

\bl Define the following cocycles in $\Hom_K(K\Gamma_1,\Lambda)$. 
\begin{align*}
\varphi_t&=(\alpha(\beta\gamma\alpha)^t,0,0,0)\quad\text{ if } 1\ppq t\ppq k-1\\
\varphi_0&=(0,\beta,0,0)\\
\theta_r&=(0,0,0,\eta^{r+1})\quad\text{ if } 0\ppq r\ppq s-1\\
\psi&=((\alpha\beta\gamma)^k,0,0,0)\\
\omega&=(\alpha+c(\beta\gamma\alpha)^{k-1}\beta\gamma+c(\alpha\beta\gamma)^k,-\beta,0,0)\end{align*}
Then a basis of $\HH^1(SD(2\B)_1^{k,s}(c))$ is given by the cohomology classes of the elements in $\B\cup\B'$ where $\B=\set{\varphi_t,\theta_r,\psi,\omega;1\ppq t\ppq k-1,1\ppq r\ppq s-1}$ and $\B'=
\begin{cases}
\set{\varphi_0,\theta_0}&\text{if $k$ and $s$ are both $0$ in $K$,}\\
\set{\varphi_0}&\text{if $k$ is $0$ and $s$ is not $0$ in $K$,}\\
\set{\theta_0}&\text{if $k$ is not $0$ and $s$ is $0$ in $K$,}\\
\varnothing&\text{if $ks$ is not $0$ in $K$.}
\end{cases}
$

The (potentially) non-zero brackets are given by
\begin{align*}
[\varphi_t,\varphi_{t'}]&=(t'-t)\varphi_{t+t'} \text{ if }0\ppq t,t'\ppq k-1,\ t+t'\ppq k-1, \\
[\theta_r,\theta_{r'}]&=(r'-r)\theta_{r+r'}  \text{ if }0\ppq r,r'\ppq s-1, \ r+r'\ppq s-1, \\
[\psi,\omega]&=\psi.
\end{align*} 
\el

It is easy to check that the Lie algebras $\HH^1(SD(2\B)_1^{k,s}(c))$ and  $\HH^1(SD(2\B)_1^{k',s'}(c'))$ are isomorphic if $\set{k,s}=\set{k',s'}$ and $c,c'\in\set{0,1}$. The Lie algebra structure does not provide any new information at this point, however is will be useful in order to distinguish the algebras $SD(2\B)_1^{k,s}(c)$ and $SD(2\B)_2^{k,s}(c)$ later.

\medskip

\paragraph{\bfseries The situation when $\car(K)\neq 2,3$ is similar} Nevertheless, we give the Lie algebra structure, since it will be used later.

\bl  Define the following cocycles in $\Hom_K(K\Gamma_1,\Lambda)$. 
\begin{align*}
\varphi_t&=(\alpha(\beta\gamma\alpha)^t,0,0,0)\quad\text{ if } 1\ppq t\ppq k-1\\
\varphi_0&=(0,\beta,0,0)\\
\theta_r&=(0,0,0,\eta^{r+1})\quad\text{ if } 0\ppq r\ppq s-1\\
\psi&=((\alpha\beta\gamma)^k,0,0,0)\\
\omega&=(ks\alpha+ksc(\beta\gamma\alpha)^{k-1}\beta\gamma,(3-k)s\beta,0,3k\eta).\end{align*}
Then a basis of $\HH^1(SD(2\B)_1^{k,s}(c))$ is given by the cohomology classes of the elements in $\B\cup\B'$ where $\B=\set{\varphi_t,\theta_r,\psi;1\ppq t\ppq k-1,1\ppq r\ppq s-1}$ and $\B'=
\begin{cases}
\set{\varphi_0,\theta_0}&\text{if $k$ and $s$ are both $0$ in $K$,}\\
\set{\omega}&\text{if $k$  and $s$ are not both $0$ in $K$.}
\end{cases}
$

The (potentially) non-zero brackets are given by
\begin{align*}
[\varphi_t,\varphi_{t'}]&=(t'-t)\varphi_{t+t'} \text{ if }0\ppq t,t'\ppq k-1,\ t+t'\ppq k-1,\\
[\theta_r,\theta_{r'}]&=(r'-r)\theta_{r+r'}  \text{ if }0\ppq r,r'\ppq s-1, \ r+r'\ppq s-1, \\
[\omega,\varphi_t]&=3ts\varphi_t \text{ for $1\ppq t\ppq k-1$, if $k\neq0$ or $s\neq0$ in $K$}\\
[\omega,\theta_r]&=3rk\theta_r \text{ for $1\ppq r\ppq s-1$, if $k\neq0$ or $s\neq0$ in $K$}\\
[\omega,\psi]&=2ks\psi \text{ (only if $k\neq0$ or $s\neq0$ in $K$)}.
\end{align*} 
\el

\subsubsection{The first Hochschild cohomology group of $SD(2\B)_2^{k,s}(c)$}

Let $\Lambda $ be the algebra $SD(2\B)_2^{k,s}(c)$ and let $\Gamma$ be the quiver of type $2\B$.

\bl\label{lemma:proj resol semidihedral2 2 simples} Define a sequence of $\Lambda $-$\Lambda $-bimodules 
\begin{align*}
Q^2=\bigoplus_{z\in Z}^2\Lambda e_{i(z)}\ot e_{t(z)}\Lambda  \stackrel{\partial^2}{\longrightarrow}&Q^1=\bigoplus_{\delta\in\Gamma_1}^2\Lambda e_{i(\delta)}\ot e_{t(\delta)}\Lambda 
\stackrel{\partial^1}{\longrightarrow}\\&\stackrel{\partial^1}{\longrightarrow}Q^0=\bigoplus_{i=1}^2\Lambda e_i\ot e_i\Lambda 
\stackrel{\partial^0}{\longrightarrow}\Lambda \rightarrow0
\end{align*} as follows. The set $Z$ is $\set{\alpha^2-c(\beta\gamma\alpha)^k,\beta\eta-\alpha\beta(\gamma\alpha\beta)^{k-1},\eta\gamma-\gamma\alpha(\beta\gamma\alpha)^{k-1},\gamma\beta-\eta^{s-1}}$ if $s>2$ and $\set{\alpha^2-c(\alpha\beta\gamma)^k,\beta\gamma\beta-(\alpha\beta\gamma)^{k-1}\alpha\beta,\gamma\beta\gamma-(\gamma\alpha\beta)^{k-1}\gamma\alpha}$ if $s=2.$ The map  $\partial ^0$ is multiplication, $\partial ^1$ is determined by $\partial ^1(e_{i(\delta)}\ot e_{t(\delta)})=\delta\ot e_{t(\delta)}-e_{i(\delta)}\ot \delta$ for $\delta\in\Gamma_1$ and $\partial ^2$ is determined by \[ 
\begin{aligned}
e_1\ot e_1\mapsto& e_1\ot \alpha+\alpha\ot e_1-c\sum_{t=0}^{k-1}(\alpha\beta\gamma)^t\left(\alpha\ot\gamma+\alpha\beta\ot e_1+e_1\ot \beta\gamma\right)(\alpha\beta\gamma)^{k-1-t}\\
e_1\ot e_2\mapsto& e_1\ot\eta+\beta\ot e_2-\sum_{t=0}^{k-1}(\alpha\beta\gamma)^t\left(e_1\ot \beta+\alpha\ot e_2\right)(\gamma\alpha\beta)^{k-1-t}\\&\qquad-\sum_{t=0}^{k-2}(\alpha\beta\gamma)^t\alpha\beta\ot \alpha\beta(\gamma\alpha\beta)^{k-2-t}\\
e_2\ot e_1\mapsto& e_2\ot\gamma+\eta\ot e_1-\sum_{t=0}^{k-1}(\gamma\alpha\beta)^t\left(e_2\ot \alpha+\gamma\ot e_1\right)(\beta\gamma\alpha)^{k-1-t}\\&\qquad-\sum_{t=0}^{k-2}(\gamma\alpha\beta)^t\gamma\alpha\ot \gamma\alpha(\beta\gamma\alpha)^{k-2-t}\\
e_2\ot e_2\mapsto& e_2\ot\beta+\gamma\ot e_2-\sum_{ r=0}^{s-2}\eta^ r\ot\eta^{s-2- r}
\end{aligned}
 \]
 if $s>2$ and by 
 \[ 
\begin{aligned}
e_1\ot e_1&\mapsto e_1\ot \alpha+\alpha\ot e_1-c\sum_{t=0}^{k-1}(\beta\gamma\alpha)^t\left(e_1\ot \gamma\alpha+\beta\ot \alpha+\beta\gamma\ot e_1\right)(\beta\gamma\alpha)^{k-1-t}\\
e_1\ot e_2&\mapsto e_1\ot\gamma\beta+\beta\ot \beta+\beta\gamma\ot e_2-\sum_{t=0}^{k-1}(\alpha\beta\gamma)^t\left(e_1\ot \beta+\alpha\ot e_2\right)(\gamma\alpha\beta)^{k-1-t}\\&\qquad-\sum_{t=0}^{k-2}(\alpha\beta\gamma)^t\alpha\beta\ot \alpha\beta(\gamma\alpha\beta)^{k-2-t}\\
e_2\ot e_1&\mapsto e_2\ot\beta\gamma+\gamma\ot \gamma+\gamma\beta\ot e_1-\sum_{t=0}^{k-1}(\gamma\alpha\beta)^t\left(e_2\ot \alpha+\gamma\ot e_1\right)(\beta\gamma\alpha)^{k-1-t}\\&\qquad-\sum_{t=0}^{k-2}(\gamma\alpha\beta)^t\gamma\alpha\ot \gamma\alpha(\beta\gamma\alpha)^{k-2-t}
\end{aligned}
 \] if $s=2.$

Then this sequence is the beginning of a minimal projective  $\Lambda $-$\Lambda $-bimodule resolution of $\Lambda .$ 
\el

\bpf The proof is the same as that of Lemma \ref{lemma:proj resol semidihedral 2 simples}.
\epf

Using the resolution above, we may now compute the dimension of $\HH^1(\Lambda )$. The proof is straightforward and is omitted.

\bp Let  $\Lambda $ be the algebra $SD(2\B)_2^{k,s}(c)$. Then 
\[ \dim\HH^1(\Lambda )=
\begin{cases}
k+s+3-c&\text{ if }\car(K)=2 \text{ and $k$ and $s$ are both even } \\
k+s+2-c&\text{ if }\car(K)=2 \text{ and $k+s$ is odd} \\
k+s+2-2c&\text{ if }\car(K)=2 \text{ and $k$ and $s$ are both odd} \\
k+s+1&\text{ if }\car(K)\neq 2 \text{ and $k$ and $s$ are  both $0$ in $K$} \\
k+s&\text{ if }\car(K)\neq 2 \text{ and  $k$ and $s$ are not both $0$ in $K$} \\
\end{cases}
 \]
\ep

\bc\label{cor:sd2 car2} Assume that $\car(K)=2$ and that $\set{k',s'}=\set{k,s}$. Then there is no stable equivalence of Morita type between $SD(2\B)^{k,s}_2(0)$ and $SD(2\B)^{k',s'}_2(1).$

If moreover $k$ or $s$ is even, then for $c\in\set{0,1}$  there is no stable equivalence of Morita type between $SD(2\B)^{k,s}_2(1)$ and $SD(2\B)^{k',s'}_1(c).$

If instead $k$ and $s$ are both odd, then for $c,c'\in\set{0,1}$ not both equal to $0,$  there is no stable equivalence of Morita type between $SD(2\B)^{k,s}_2(c)$ and $SD(2\B)^{k',s'}_1(c').$

Finally, if $\car(K)=3$,  $\set{k',s'}=\set{k,s}$ with $ks=0$ in $K$ and if $c,c'\in\set{0,1}$, then  there is no stable equivalence of Morita type between $SD(2\B)^{k,s}_2(c)$ and $SD(2\B)^{k',s'}_1(c').$
\ec

In order to go further, we now consider the Lie algebra structure of $\HH^1(\Lambda )$.  Once more, we identify a morphism $f\in\Hom_{\Lambda -\Lambda }(Q^1,\Lambda )$ with  $g\in\Hom_{K\Gamma_0-K\Gamma_0}(K\Gamma_1,\Lambda )$ and with the quadruple $(g(\alpha),g(\beta),g(\gamma),g(\eta)).$

\paragraph{\bfseries First assume that $\car(K)=2$} We start with a basis for  $\HH^1(SD(2\B)^{k,s}_2(c))$. 

\bl Assume that $s>2$. Define the following cochains in $\Hom_K(K\Gamma_1,\Lambda)$.
\begin{align*}
&\varphi_t=(\alpha(\beta\gamma\alpha)^t,0,0,0)\quad\text{ if } 1\ppq t\ppq k-1,&\qquad&\varphi_0=(\alpha,0,0,0),
\\
&\theta_ r =(0,0,0,\eta^{ r +1})\quad\text{ if } 2\ppq r \ppq s-1,&&\psi_1=((\alpha\beta\gamma)^k,0,0,0),\\
&\theta_1=(0,(s-1)\alpha\beta(\gamma\alpha\beta)^{k-1},0,\eta^2),&&\psi_0=(e_1,0,0,(\gamma\alpha\beta)^{k-1}),
\\
&\omega=((\beta\gamma\alpha)^{k-1}\beta\gamma,0,0,0),&&\theta_0=(\alpha,\beta,0,\eta),
\\
&\chi=(\alpha,0,0,\eta).
\end{align*}
Set \begin{align*}
\B&=\set{\varphi_t;\theta_r;\psi_1;\omega;1\ppq t\ppq k-1,1\ppq r\ppq s-1}\\&\B'=
\begin{cases}
\set{\varphi_0,\theta_0}&\text{ if $k$ and $s$ are even}\\
\set{\theta_0}&\text{if $k$ is even and $s$ is odd}\\
\set{\varphi_0}&\text{if $s$ is even and $k$ is odd}\\
\set{\chi}&\text{if $ks$ is odd and $c=0$}\\
\varnothing&\text{if $ks$ is odd and $c=1$},
\end{cases}\\
\text{and }&\B''=
\begin{cases}
\set{\psi_0}&\text{if $c=0$}\\
\varnothing&\text{if $c=1$}.
\end{cases}
\end{align*} Then the cohomology classes of the elements in  $\B\cup \B'\cup \B''$ form a basis for $\HH^1(SD(2\B)_2^{k,s}(c)).$
\el

\bl When $s=2$, a basis is given by the non-zero cohomology classes of the following elements of $\Hom_K(K\Gamma_1,\Lambda)$, written as $(g(\alpha),g(\beta),g(\gamma))$:
\begin{align*}
&\varphi_t=(\alpha(\beta\gamma\alpha)^t,0,0,0)\quad\text{ if } 1\ppq t\ppq k-1,&\qquad&\theta_1=(0,(\alpha\beta\gamma)^{k-1}\alpha\beta,0),\\
&\psi_1=((\alpha\beta\gamma)^k,0,0),&&\omega=((\beta\gamma\alpha)^{k-1}\beta\gamma,0,0),\\
&\psi_0=(1-c)(e_1,(\alpha\beta\gamma)^{k-2}\alpha\beta,0),&&\\
&\varphi_0=(\alpha,k\beta,0),&&\theta_0=(1-k)(0,\beta,0).
\end{align*} 
\el

We may now compute the brackets of these basis elements.

\bl\label{lemma:brackets sd2 car 2} We use the notation in the previous lemmas. We describe the (potentially) non-zero brackets in $\HH^1(SD(2\B)^{k,s}_2(c))$. 

For all $k $ and $s$ we have
\begin{align*}
[\varphi_t,\varphi_{t'}]&=(t+t')\varphi_{t+t'} \ (t,t'\pgq 1),&[\theta_r,\theta_{r'}]&=(r+r')\theta_{r+r'}  \ (r,r'\pgq 1),\\
[\varphi_t,\varphi_0]&=t\varphi_t,&[\theta_r,\theta_0]&=r\theta_r,\\
[\varphi_1,\psi_0]&=(1-c)(1-k)\theta_{s-1},&[\theta_1,\psi_0]&=(1-c)(1-s)\varphi_{k-1},\\
[\psi_1,\varphi_0]&=\psi_1,&[\psi_1,\theta_0]&=\psi_1,\\
[\psi_0,\varphi_0]&=(1-c)\psi_0,&[\psi_0,\theta_0]&=(1-c)\psi_0,\\
[\psi_1,\psi_0]&=(1-c)\omega,&&
\end{align*} and, if $ks\neq 0$ in $K$ and $c=0$, 
\begin{align*}
[\varphi_t,\chi]&=t\varphi_t\ (t\pgq 1),&[\theta_r,\chi]&=r\theta_r\ (r\pgq 1),\\
[\omega,\chi]&=\omega,&[\psi_0,\chi]&=\psi_0.
\end{align*}
\el

\br It is easy to check that $\HH^1(SD(2\B)_2^{k,s}(c))$ and  $\HH^1(SD(2\B)_2^{s,k}(c))$ are isomorphic Lie algebras. 
\er

\bc\label{cor:sd21 car2} Assume that $\car(K)=2$. Let $k,k',s,s'$ be integers such that $\set{k,s}=\set{k',s'}$ and let $c$, $c'$ be in $\set{0,1}$. Suppose that one of the following holds:
\begin{enumerate}[(i)]
\item $k$ and $s$ are both even and $cc'=0$;
\item $k+s$ is odd;
\item $ks$ is odd and $(c,c')\neq (0,0)$.
\end{enumerate}  Then  there is no stable equivalence of Morita type between the algebras  $SD(2\B)_2^{k,s}(c)$ and $SD(2\B)_1^{k',s'}(c')$. 

\ec

\bpf Set  $\mathfrak{g} =\HH^1(SD(2\B)_2^{k,s}(c))$ and $\mathfrak{g} '=\HH^1(SD(2\B)_1^{k',s'}(c'))$. Let $\L^i(\mg)$ and $\L^i(\mg')$ be the $i$\tup{th} term in the lower central series of $\mathfrak{g}$ and $\mathfrak{g}'$ respectively. Write the basis elements in $\mathfrak{g}'$ with dashes.

In case (i) $\L^2(\mg)$ is the span of the set $\set{\varphi_{2p+1};\theta_{2q+1};\psi_1;(1-c)\omega;(1-c)\psi_0;0\ppq p\ppq \frac{k}{2}-1,0\ppq q\ppq \frac{s}{2}-1}$ and its dimension is $\frac{k+s}{2}+3-2c$, and $\L^2({\mathfrak{g}'})$ is the span of the set $\set{\varphi'_{2p+1};\theta'_{2q+1};c'\omega';0\ppq p\ppq \frac{k'}{2}-1,0\ppq q\ppq \frac{s'}{2}-1}$ and it has dimension $\frac{k+s}{2}+c'$. These dimensions are different when $cc'=0$, therefore $\mathfrak{g}$ and $\mathfrak{g}'$ are not isomorphic.

In case (ii), we may assume that $k=k'$ is odd and $s=s'$ is even. Here,  $\L^1(\mg)$ is the span of $\set{\varphi_{2p+1};\theta_{2q+1};\psi_1;(1-c)\varphi_{k-1};(1-c)\omega;(1-c)\psi_0;1\ppq p\ppq \frac{k-3}{2},0\ppq q\ppq \frac{s}{2}-1}$ so its dimension is $\frac{k+s-1}{2}+3-3c+\delta_{k,1}$, and $\L^1(\mg')$ is the span of $\set{\varphi'_{2p+1};\theta'_{2q+1};\omega';1\ppq p\ppq \frac{k-3}{2},0\ppq q\ppq \frac{s}{2}-1}$ and has dimension $\frac{k+s-1}{2}+\delta_{k,1}$, and these dimensions are different when $(c,c')\neq (1,0)$. Moreover, if $i>\frac{k-3}{2}$, we have $\dim \L^i(\mg)=\frac{s}{2}+4-3c$ and $\dim \L^i(\mg')=\frac{s}{2}+c'$, which are different when $c=1$ and $c'=0$.  Therefore $\mathfrak{g}$ and $\mathfrak{g}'$ are not isomorphic.

\sloppy Finally, in case (iii), again assume that $k=k'$ and $s=s'$ for the proof. In this case,  $\L^1(\mg)$ is the span of $\set{\varphi_{2p+1};\theta_{2q+1};(1-c)\omega;(1-c)\psi_0;(1-c)\varphi_1;(1-c)\theta_1;1\ppq p\ppq \frac{k-3}{2},1\ppq q\ppq \frac{s-3}{2}}$ and $\L^1(\mg')$ is spanned by \linebreak $\set{\varphi'_{2p+1};\theta'_{2q+1};\omega';(1-c')\chi';(1-c)\varphi'_1;(1-c')\theta'_1;1\ppq p\ppq \frac{k-3}{2},1\ppq q\ppq \frac{s-3}{2}}$, therefore $\dim \L^1(\mg)-\dim \L^1(\mg')=(c'-c)\delta_{k,1}+3c'-4c$, which is non-zero when $(c,c')\neq (0,0)$. Therefore  $\mathfrak{g}$ and $\mathfrak{g}'$ are not isomorphic.
\epf

\br If $ks$ is odd, then the Lie algebras $\HH^1(SD(2\B)_1^{k,s}(0))$ and  $\HH^1(SD(2\B)_2^{k,s}(0))$ are isomorphic, so that the Lie algebra structure of the first Hochschild cohomology group does not bring anything new. Indeed, if $\set{\varphi_t;\theta_ r;\psi;\chi;\omega;\zeta_0;1\ppq t\ppq k-1;1\ppq r\ppq s-1}$ is our basis of $\HH^1(SD(2\B)_1^{k,s}(0))$  and $\set{\varphi_t';\theta_r';\psi_0';\psi_1';\omega';\chi';1\ppq t\ppq k-1,1\ppq r\ppq s-1}$ is our basis of  $\HH^1(SD(2\B)_2^{k,s}(0))$, then the isomorphism is determined by 
\begin{align*}
&\varphi_t\mapsto \varphi_t',&&\theta_r\mapsto \theta_r',&&\omega\mapsto \omega',\\
&\psi\mapsto \psi'_1,&&\zeta_0\mapsto\chi',&&\chi\mapsto\psi_0'.
\end{align*}
\er

\paragraph{\bfseries We now assume that $\car(K)\neq 2$}
In order to differentiate  the algebras of type $SD(2\B)_1$ and $SD(2\B)_2$ up to stable equivalence of Morita type,  we give the Lie algebra structure when $\car(K)\neq 2$.

\bl\label{lemma:basis-brackets sd2 car not 2}  We define cocycles in $\Hom_{K\Gamma_0-K\Gamma_0}(K\Gamma_1,\Lambda )$ as follows.
\begin{align*}
\varphi_t&=(\alpha(\beta\gamma\alpha)^t,0,0,0)\quad\text{ if } 1\ppq t\ppq k-1\\
\theta_r&=(0,0,0,\eta^{r+1})\quad\text{ if } 2\ppq r\ppq s-1\\
\theta_1&=(0,(s-1)(\alpha\beta\gamma)^{k-1}\alpha\beta,0,\eta^2)\\
\psi&=((\alpha\beta\gamma)^k,0,0,0)
\\
\varphi_0&=(\alpha-c(\beta\gamma\alpha)^{k-1}\beta\gamma,0,0,0)\text{ and } \theta_0=(0,\beta,0,\eta)-\varphi_0\text{ if $k=0$ and $s=0$ in $K$,}\\
\omega&=(2(k+s-ks)\alpha+c(3ks-2k-2s)(\beta\gamma\alpha)^{k-1}\beta\gamma,2k(s-1)\beta,0,2k\eta)\\&\qquad\qquad\text{ if $k\neq 0$ or $s\neq0$ in $K$.}
\end{align*}
 Set $\B=\set{\varphi_t;\theta_ r;\psi;1\ppq t\ppq k-1;1\ppq r\ppq s-1}$ and $\B'=
\begin{cases}
\set{\varphi_0;\theta_0}&\text{ if $k$ and $s$ are both zero in $K$,}\\
\set{\omega}&\text{ if $k\neq 0$ or $s\neq 0$ in $K$.}
\end{cases}
$ 

Then the cohomology classes of the elements in $\B\cup\B'$ form a basis of $\HH^1(SD(2\B)_2^{k,s}(c))$.

The (potentially) non-zero brackets are given by
\begin{align*}
[\varphi_t,\varphi_{t'}]&=(t'-t)\varphi_{t+t'}\rlap{ $\text{ if }1\ppq t,t'\ppq k-1,\ t+t'\ppq k-1,\ $}\\
[\theta_r,\theta_{r'}]&=(r'-r)\theta_{r+r'} \rlap{ $ \text{ if }1\ppq r,r'\ppq s-1, \ r+r'\ppq s-1,\ $}\\
[\varphi_0,\varphi_t]&=t\varphi_t,&[\theta_0,\theta_r]&=r\theta_r,\\
[\psi,\varphi_0]&=\psi,&[\psi,\theta_0]&=-\psi,\\
[\omega,\varphi_t]&=2st\varphi_t,&[\omega,\theta_r]&=2kr\theta_r,\\
[\omega,\psi]&=2(2ks-k-s)\psi.
\end{align*} 
\el

\br Here again, if we specialise to  $s=2$, a basis of $\HH^1(SD(2\B)_2^{k,2}(c))$ is given by the non-zero cohomology classes of the following elements of $\Hom_K(K\Gamma_1,\Lambda)$, written as $(g(\alpha),g(\beta),g(\gamma))$:
\begin{align*}
&\varphi_t=(\alpha(\beta\gamma\alpha)^t,0,0,0)\quad\text{ if } 1\ppq t\ppq k-1,&\qquad&\theta_1=(0,(\alpha\beta\gamma)^{k-1}\alpha\beta,0),\\
&\omega=((2-k)\alpha+2c(k-1)(\beta\gamma\alpha)^{k-1}\beta\gamma,k\beta,0),&&\psi=((\alpha\beta\gamma)^k,0,0)
\end{align*} and the brackets are the same as those given in Lemma \ref{lemma:basis-brackets sd2 car not 2} above.
\er

It is easy to check that the Lie algebras $\HH^1(SD(2\B)_2^{k,s}(c))$ and  $\HH^1(SD(2\B)_2^{k',s'}(c'))$ are isomorphic if $\set{k,s}=\set{k',s'}$ and $c,c'\in\set{0,1}$. The Lie algebra structure does not provide any new information within the family $SD(2\B)_2^{k,s}(c)$, but we have the following result.

\bc\label{cor:sd21 car not 2} Assume that $\car(K)=3$, that $\set{k,s}=\set{k',s'}$ and that $c,c'\in\set{0,1}.$ Then  there is no stable equivalence of Morita type between $SD(2\B)_2^{k,s}(c)$ and $SD(2\B)_1^{k',s'}(c')$. 

Assume that $\car(K)\neq 2,3$ and that either $ks=0$ in $K$, or that $ks\neq 0$ and $2ks-k-s=0$ in $K$. Then  there is no stable equivalence of Morita type between $SD(2\B)_2^{k,s}(c)$ and $SD(2\B)_1^{k',s'}(c')$. 
\ec

\bpf  Set  $\mg =\HH^1(SD(2\B)_2^{k,s}(c))$ and $\mg '=\HH^1(SD(2\B)_1^{k',s'}(c'))$. Let $\L^i(\mg)$ (respectively $\L^i(\mg')$) denote  the $i$\tup{th} term in the lower central series of $\mg$ (respectively $\mg'$). 
\begin{itemize}
\item First assume that $\car(K)=3$. We already know from Corollary \ref{cor:sd2 car2} that  there is no stable equivalence of Morita type between $SD(2\B)_2^{k,s}(c)$ and $SD(2\B)_1^{k',s'}(c')$ when $ks=0$ in $K$.

Therefore, assume that  that $ks\neq 0$ in $K$. 
Then the centre of the Lie algebra $\mg$ is spanned by $\psi$ if $(k,s)=(1,1)$ in $K^2$ and vanishes otherwise, so its dimension is at most $1, $ whereas the centre  of the Lie algebra $\mg'$ is spanned by $\set{\varphi_{k-1},\theta_{s-1}}$ so its dimension is at least $2.$
Therefore the algebras $\mg$ and $\mg'$ are not isomorphic and the first part of the corollary follows.

\item If $\car(K)\neq 2,3$ and  $ks=0$, then $\dim \L^1_\mg=\dim \L^1_{\mg'}+1$ (the extra element is $\psi$), hence the Lie algebras $\mg$ and $\mg'$ are not isomorphic.

\item  If $\car(K)\neq 2,3$ and   $ks\neq 0$ and $k+s-2ks=0$ in $K$, then the centre of $\mg'$ is zero, whereas that of $\mg$ is spanned by $\psi$ and has dimension $1$, hence the Lie algebras $\mg$ and $\mg'$ are not isomorphic.
\qedhere
\end{itemize}
\epf

\bp\label{prop:sd21 constants non zero}
Assume that  $\car(K)\neq 2,3$ and that $\set{k,s}=\set{k',s'}$. Put $\lambda=3^{-1}2ks\neq 0$ and $\mu=2ks-k-s$ and assume that $\lambda\mu\neq 0$ in $K$ and that the following subsets of $K$ are not equal:
\begin{align*}
&\mathscr{E}_\lambda=\set{s\lambda^{-1},2s\lambda^{-1},k\lambda^{-1},2k\lambda^{-1},(s\lambda^{-1})^{-1},(2s\lambda^{-1})^{-1},(k\lambda^{-1})^{-1},(2k\lambda^{-1})^{-1}}\\
\text{and } &\mathscr{E}_\mu=\set{s\mu^{-1},2s\mu^{-1},k\mu^{-1},2k\mu^{-1},(s\mu^{-1})^{-1},(2s\mu^{-1})^{-1},(k\mu^{-1})^{-1},(2k\mu^{-1})^{-1}}.
\end{align*} Then  there is no stable equivalence of Morita type between $SD(2\B)_2^{k,s}(c)$ and $SD(2\B)_1^{k',s'}(c')$. 
\ep

\bpf
For $\lambda\in K^*$, let $\mg_\lambda$ be the $6$-dimensional Lie algebra with basis $\set{e_0,\ldots,e_5}$ and whose bracket is determined by $[e_0,e_i]=\nu_ie_i$ with $(\nu_1,\nu_2,\nu_3,\nu_4,\nu_5)=(s,2s,k,2k,\lambda)$.

Now consider the algebra $\HH^1(SD(2\B)_2^{k,s}(c))/\D^2(\HH^1(SD(2\B)_2^{k,s}(c)))$ where $\D^2(\mg)$ is the second term in the derived series of $\mg$. Since $\D^2(\HH^1(SD(2\B)_2^{k,s}(c)))$ is spanned by the $\varphi_t$ and $\theta_r$ for $t\pgq 3$ and $r\pgq 3$,  this is a Lie algebra that is isomorphic to $\mg_\lambda$ with $\lambda=2ks-k-s$, via the isomorphism given by $e_0=\frac{1}{2}\omega$, $e_1=\varphi_1$, $e_2=\varphi_2$, $e_3=\theta_1$, $e_4=\theta_2$ and $e_5=\psi$  (recall that $2$ and $3$ are invertible in $K$). Similarly, the Lie algebra $\HH^1(SD(2\B)_1^{k,s}(c'))/\D^2(\HH^1(SD(2\B)_1^{k,s}(c')))$ is isomorphic to $\mg_\mu$ with $\mu=\frac{2ks}{3}$ (the isomorphism sends $\omega$ to $\frac{1}{3}\omega$ in this case). 

If $\HH^1(SD(2\B)_2^{k,s}(c))$ and $\HH^1(SD(2\B)_1^{k,s}(c'))$ are isomorphic, then so are $\mg_\lambda$ and $\mg_\mu$. 

We now prove that if the sets $\mathscr{E}_\lambda$ and $\mathscr{E}_\mu$ are distinct, then the Lie algebras $\mg_\lambda$ and $\mg_\mu$ are not isomorphic,  using generalised derivations.

For $\rho\in K^*$, we consider $\der_{\mg_\lambda}(\rho,1,1).$ Let $D$ be a $(\rho,1,1)$-derivation of $\mg_\lambda$. Set $D(e_j)=\sum_{i=0}^5a_{ij}e_i$ for $i=0,1,\ldots,5$. Then, for $i=1,2,\ldots,5$, we have 
\[ \rho D([e_0,e_i])=[D(e_0),e_i]+[e_0,D(e_i)], \]
 which is equivalent to the set of equations 
\[ 
\begin{cases}
\rho\nu_ia_{0i}=0&\text{ for }1\ppq i\ppq 5,\\
(\rho\nu_i-\nu_j)a_{ji}=0&\text{ for }1\ppq i\neq j\ppq 5,\\
(\rho-1)a_{ii}=a_{00}&\text{ for }1\ppq i\ppq 5
\end{cases} 
 \] that is equivalent to \[ 
\begin{cases}
a_{0i}=0&\text{ for }1\ppq i\ppq 5,\\
(\rho\nu_i-\nu_j)a_{ji}=0&\text{ for }1\ppq i\neq j\ppq 5,\\
(\rho-1)a_{ii}=a_{00}&\text{ for }1\ppq i\ppq 5.
\end{cases} 
 \] Note that the equations that come from the identities $\rho D([e_j,e_i])=[D(e_j),e_i]+[e_j,D(e_i)]$ for $1\ppq i\neq j\ppq 5$ are a consequence of the first five equations above. Therefore these equations characterise $D$. The subset of the equations above that involve the parameter $\lambda$ is 
\begin{align*}
&(\rho\lambda-s)a_{15}=0&&(\rho s-\lambda)a_{51}=0\\
&(\rho\lambda-2s)a_{25}=0&&(2\rho s-\lambda)a_{52}=0\\
&(\rho\lambda-k)a_{35}=0&&(\rho k-\lambda)a_{53}=0\\
&(\rho\lambda-2k)a_{45}=0&&(2\rho k-\lambda)a_{54}=0.
\end{align*} 
Therefore, if $\mu\in K^*$ is another parameter, and if $\rho\in\mathscr{E}_\mu$ and  $\rho\not\in\mathscr{E}_\lambda$, then there are strictly fewer equations characterising $\der_{\mg_\mu}(\rho,1,1)$ than those  characterising $\der_{\mg_\lambda}(\rho,1,1)$. It follows that $\dim \der_{\mg_\lambda}(\rho,1,1)<\dim \der_{\mg_\mu}(\rho,1,1)$ and hence that $\mg_\lambda$ and $\mg_\mu$ are not isomorphic Lie algebras.
\epf

Finally, Theorem \ref{thm:sd} is obtained by combining Propositions \ref{prop:sd1 car2} and \ref{prop:sd21 constants non zero} and  Corollaries \ref{cor:sd1 car2}, \ref{cor:sd2 car2}, \ref{cor:sd21 car2} and \ref{cor:sd21 car not 2}.

\section{Algebras of quaternion type}\label{sec:quaternion}

As we mentioned in Subsection \ref{subsec:questions:quaternion}, we shall only consider the  local tame symmetric algebras of quaternion
type. Using a result of  Erdmann and  Skowro\'nski,   in this case we can compute  the dimensions of all the Hochschild cohomology
groups. The dimension of the first Hochschild cohomology group, as well as the Lie algebra structure of the first cohomology group $\HH^1(\Lambda)$,
give new information on stable equivalence of Morita type, but we are not able to distinguish all
the algebras. The main result of this subsection is Corollary \ref{cor:classification q 2}.

Once more, we assume that the field $K$ has characteristic $2$. We have defined the algebras $Q(1\A)_1^k$ and $Q(1\A)_2^k(c,d)$ in Subsection \ref{subsec:questions:list algebras}. In these algebras, the following relations hold: $x^3=(xy)^k=(yx)^k=y^3$ and $x^4=0.$  The element $z:=(xy)^{k-1}+(yx)^{k-1}$ is central in these algebras by \cite{E}, therefore from the equalities $y^2z=yzy=zy^2$, using the other relations, we obtain $x^2y=0=yx^2$. It then follows that $xy^2=0=y^2x$, and that $y(xy)^k=y^4=0$, even  in $Q(1\A)_1^k$. We may therefore view  $Q(1\A)_1^k$ as $Q(1\A)_2^k(0,0)$.

\sloppy Fix an integer $k\pgq2.$ For $(c,d)\in K^2$, let $I^k(c,d)$ be the ideal in $K\rep{x,y}$ generated by
 the set  $\set{(xy)^k+(yx)^k;x^2+(yx)^{k-1}y+c(xy)^k;y^2+(xy)^{k-1}x+d(yx)^k;(xy)^kx;(yx)^ky}$ and let $\Lambda:=K\rep{x,y}/I^k(c,d)$ be one of the algebras  $Q(1\A)_1^k$ or $Q(1\A)_2^k(c,d)$. 
Clearly, $Q(1\A)_2^k(c,d)\cong Q(1\A)_2^k(d,c)$.

\subsection{Dimensions of the Hochschild cohomology groups}

  Erdmann and
 Skowro\'nski have shown in \cite{ES} that $\Lambda$ is  periodic of period $4$ and they give explicitly  a minimal projective resolution
of $\Lambda$ as a $\Lambda$-$\Lambda$-bimodule in \cite[Theorem 5.9]{ES}: 
\begin{equation}\label{eq:ES resol local}
\begin{aligned}
 0\rightarrow \Lambda \stackrel j\rightarrow \Lambda\ot \Lambda\stackrel{{\partial ^3}}{\rightarrow}( \Lambda\ot \Lambda)^2&=(\Lambda\s{x}\Lambda)\oplus(\Lambda\s{y}\Lambda)\stackrel{{\partial ^2}}{\rightarrow}\\&\stackrel{{\partial ^2}}{\rightarrow}(\Lambda\s{x}\Lambda)\oplus(\Lambda\s{y}\Lambda )\stackrel{{\partial ^1}}{\rightarrow}\Lambda\ot \Lambda \stackrel{\partial ^0}{\rightarrow}\Lambda\rightarrow0
\end{aligned}\end{equation} where $\partial ^0$ is multiplication, ${\partial ^1}(1\s{\alpha}1)=\alpha\ot 1+1\ot \alpha$ for $\alpha\in\set{x,y}$ and ${\partial ^2}$, ${\partial ^3}$ and $j$ are determined by: 
\begin{align*}
&{\partial ^2}(1\s{x}1)=x\s{x}1+1\s{x}x+\sum_{t=0}^{k-2}(yx)^ty\s{x}y(xy)^{k-2-t}\\&\qquad+\sum_{t=0}^{k-1}\left((yx)^t\s{y}(xy)^{k-1-t}+c(yx)^ ty\s{x}(yx)^{k-1-t}+c(yx)^t\s{y}x(yx)^{k-1-t} \right)\\
&{\partial ^2}(1\s{y}1)=y\s{y}1+1\s{y}y+\sum_{t=0}^{k-2}(xy)^tx\s{y}x(yx)^{k-2-t}\\&\qquad+\sum_{t=0}^{k-1}\left((xy)^t\s{x}(yx)^{k-1-t}+d(xy)^tx\s{y}(xy)^{k-1-t}+d(xy)^t\s{x}y(xy)^{k-1-t} \right)\\
&{\partial ^3}(1\ot1)=(x\s{x}1+1\s{x}x)(1+cx+c^2x^2)+(y\s{y}1+1\s{y}y)(1+dy+d^2y^2)\\
&j(1)=\sum_{t=0}^{k-1}\left((xy)^t\ot(xy)^{k-t}+(yx)^{t+1}\ot(yx)^{k-t-1}\right.\\&\qquad\qquad\left.+(xy)^tx\ot y(xy)^{k-1-t}+(yx)^ty\ot x(yx)^{k-1-t}\right).
\end{align*}

Again, it is straightforward to check that
\[ \dim\HH^1(\Lambda)=
\begin{cases}
k+5&\text{ if $k$ is even of if $k$ is odd and $(c,d)=(0,0)$}\\
k+4&\text{ if $k$ is odd and $(c,d)\neq(0,0)$.}
\end{cases}
 \]

We shall now give the dimensions of all the Hochschild cohomology groups for $\Lambda.$ Note that
for $\Lambda=Q(1\A)_1^k$, these were already given in \cite{G07}.

\bp \label{prop:HHn quaternion 1 simple}
We have the following dimensions
\begin{align*}
&\dim\HH^i(Q(1\A)_1^k)=
\begin{cases}
k+3 &\text{ if $i\equiv 0,3\pmod{4}$}\\
k+5 &\text{ if $i\equiv 1,2\pmod{4}$}
\end{cases}\\
&\dim\HH^i(Q(1\A)_2^k(c,d))=
\begin{cases}
k+3 &\text{ if $i\equiv 0,3\pmod{4}$}\\
k+5 &\text{ if $i\equiv 1,2\pmod{4}$ and $k$ is even}\\
k+4 &\text{ if $i\equiv 1,2\pmod{4}$ and $k$ is odd}
\end{cases}
\end{align*}
\ep

\bpf Let $\Lambda$ be one of the algebras $Q(1\A)_1^k$ or
$Q(1\A)_2^k(c,d)$. By \cite[Theorem 5.9]{ES}, we have $\Omega_{\Lambda^e}^4(\Lambda)\cong \Lambda$. In
particular, $\HH^{i+4}(\Lambda)=\HH^i(\Lambda)$ for all $i\pgq1.$ Moreover, $\Lambda$ is periodic
Frobenius of period $\pi\ppq 4$ and dimension $\pi-1$ in the sense of \cite{EuSc}. We also deduce that
$\Omega_\Lambda^\pi\cong \id_{\underline{\mathrm{mod}}\Lambda}$, so that if we assume $\pi\ppq 3$, then the
stable Calabi-Yau dimension of $\Lambda$ in the sense of \cite{ES} is at most $2.$ However, by
\cite[Proposition 5.8 and Corollary 5.10]{ES}, this last stable dimension is equal to $3.$ Therefore
$\pi=4.$

It now follows from \cite[Theorem 2.3.27(ii)]{EuSc}, using the fact that $\Lambda$ is symmetric
(hence the $K$-dual $\Lambda^*$ is isomorphic to $\Lambda$ as a $\Lambda$-$\Lambda$-bimodule) and
using Corollary 2.1.13 and Definitions 2.1.22 to 2.1.28 in \cite{EuSc} as well as the two-sided resolution of
$\Lambda$ obtained from \cite{ES}, that $\dim \sHH^{3-i}(\Lambda)=\dim \sHH^i(\Lambda)$ for
$i=0,1,2,3$ and therefore that $\dim \HH^{2}(\Lambda)=\dim \HH^1(\Lambda)$. (A direct computation using the resolution in \cite{ES} also gives 
this last fact.) Moreover, computing the dimensions from the complex obtained from \eqref{eq:ES resol local}, we get $\dim \HH^3(\Lambda)=\dim\HH^4(\Lambda)=\dim\HH^0(\Lambda)-\dim(\im(?\circ j\circ \partial ^0))=\dim\HH^0(\Lambda).$  The result follows.
\epf

We can therefore resolve some of the classification questions in this case (note that the first Hochschild cohomology group is enough for this).

\bc\label{cor:HHn quaternion 1 simple} If $k$ is odd then there is no stable equivalence of Morita type between $Q(1\A)_1^k$ and
$Q(1\A)_2^k(c,d)$.
\ec

\subsection{Lie algebra structure on $\HH^1(\Lambda)$}\label{subsection:Lie quaternion}

Let $\Gamma$ be a quiver of type $1\A$,  with arrows $x$ and $y$. Then the local tame symmetric algebras of quaternion type
 may be defined as $K\Gamma/I_2^k(c,d)$ for $(c,d)\in K^2.$ Let $\Lambda$ be such an algebra. 

We then have $\Hom_{\Lambda-\Lambda}((\Lambda\s{x}\Lambda)\oplus(\Lambda\s{y}\Lambda),\Lambda)\cong\Hom_K(K\Gamma_1,\Lambda)$ via the correspondence $f\leftrightarrow g$ given by $f(1\s{x}1)=g(x)$ and $f(1\s{y}1)=g(y).$

Moreover, if $Z=\set{x^2+(yx)^{k-1}+c(yx)^k,y^2+(xy)^{k-1}x+d(xy)^k}$, we can identify $Q^2$ with $\Lambda\ot_{K\Gamma_0}KZ\ot_{K\Gamma_0}\Lambda$.

Define the following elements in $\Hom_K(K\Gamma_1,\Lambda)$ (as pairs $(g(x),g(y))$):
\begin{align*}
&\varphi_t=
(x(yx)^t, 0)
\quad\text{ for $1\ppq t\ppq k-1$},&\quad&
\theta_1=(y(xy)^{k-1}, 0),&\quad&\theta_{-1}( 0,{x(yx)^{k-1}}),\\
&\chi=( 1+cx, x(yx)^{k-2}+d(xy)^{k-1}),&&\theta_{-2}=( 0, (xy)^k),&\quad&\theta_2( (xy)^k, 0),
\\& \omega=( y(xy)^{k-2}+c(yx)^{k-1}, 1+dy).
\end{align*}

 We then have the following result.

\bl\label{lemma:basis brackets q1} We keep the notation above.
\begin{enumerate}[(1)]
\item If $k$ is odd and $(c,d)\neq (0,0)$ then a basis for $\HH^1(Q(1\A)_2^k(c,d))$ is given by the
  cohomology classes of
  \[\set{\theta_1;\theta_{-1};\theta_2;\theta_{-2};\varphi_t,1\ppq t\ppq k-1; \psi:=d\chi+c\omega}.\]

Otherwise, $\set{\theta_1;\theta_{-1};\theta_2;\theta_{-2};\varphi_t,1\ppq t\ppq k-1;\chi;\omega}$ is a basis of $\HH^1(Q(1\A)_2^k(c,d))$  and of $\HH^1(Q(1\A)_1^k)$.

\item The (potentially) non-zero brackets of these basis elements are the following:
\begin{alignat*}{2}
& [\varphi_t,\chi]=tc\varphi_t\text{ for $t>1$,}&\quad& [\varphi_t,\omega]=td\varphi_t\text{ for $t>1$,}\\
&[\varphi_1,\chi]=c\varphi_1+(k-1)(\theta_{-1}+d\theta_{-2}),&&[\varphi_1,\omega]=d\varphi_1+(k-1)(\theta_1+c\theta_2),\\&[\varphi_t,\varphi_{t'}]=(t+t')\varphi_{t+t'},&&[\varphi_1,\psi]=c(\theta_1+c\theta_2)+d(\theta_{-1}+d\theta_{-2}),
\\
&[\theta_{1},\chi]=kc\theta_1,&\quad& [\theta_1,\omega]=kd\theta_1,\\
& [\theta_{-1},\chi]=kc\theta_{-1},&\quad&[\theta_{-1},\omega]=kd\theta_{-1},\\
& [\theta_2,\chi]=\theta_1+(k-1)c\theta_2,&\quad& [\theta_2,\omega]=\varphi_{k-1}+kd\theta_2,\\
&[\theta_{-2},\chi]=\varphi_{k-1}+kc\theta_{-2},&\quad&[\theta_{-2},\omega]=\theta_{-1}+(k-1)d\theta_{-2},\\
& [\theta_2,\psi]=c\varphi_{k-1}+d(\theta_1+c\theta_2),&\quad&[\theta_{-2},\psi]=d\varphi_{k-1}+c(\theta_{-1}+d\theta_{-2}).
\end{alignat*}\end{enumerate}
\el

We start with a special case.

\bl\label{lemma:q 2 when k=2} If $cd\neq 0$ in $K$, then for any $d'\in K$ there is no stable equivalence of Morita type between the
 algebras $Q(1\A)_2^2(0,d')$ and $Q(1\A)_2^2(c,d)$.
\el

\bpf In the basis described in Lemma \ref{lemma:basis brackets q1}, the Killing form of the Lie algebra $\HH^1(Q(1\A)_2^2(c,d))$ has matrix $\begin{pmatrix}0_5&0&0\\0&0&cd\\0&cd&0 \end{pmatrix}$. Therefore its rank is $2$ if $cd\neq 0$ and $0$ if $cd=0$. The result follows, since the rank of the Killing form invariant under an isomorphism of Lie algebras.
\epf

We then compute the first two terms in the lower central series. In view of Lemma \ref{lemma:q 2 when k=2}, we need only consider the cases where $cd=0$, that is, $(c,d)=(0,0)$ and $c=0,$ $d\neq 0$.

\bp\label{prop:lc series q1} We keep the notation above.  Then  $\L^1(\HH^1(\Lambda))$ is spanned by:
\begin{enumerate}[(a)]
\item $\begin{aligned}[t]\biggl\{\textstyle\varphi_{2p+1},1\ppq p\ppq \frac{k-3}{2};&c\varphi_{k-1}+d(\theta_1+c\theta_2);d\varphi_{k-1}+c(\theta_{-1}+d\theta_{-2});\\&c(\theta_1+c\theta_2)+d(\theta_{-1}+d\theta_{-2})\biggr\}\end{aligned}$\linebreak
  if $k$ is odd and $cd=0,$ $(c,d)\neq (0,0);$ the dimension is $\frac{k+3}{2}$;
\item $\set{\varphi_{2p+1},1\ppq p\ppq \frac{k-3}{2};\varphi_{k-1};\theta_1;\theta_{-1}}$ if $k$ is
  odd and $(c,d)=(0,0)$;  the dimension is $\frac{k+3}{2}$;

\item $\set{\varphi_{2p+1},0\ppq p\ppq \frac{k-4}{2};\varphi_{k-1};\theta_1+c\theta_2;\theta_{-1}+d\theta_{-2}}$ if $k$ is even and $cd=0,$ $(c,d)\neq (0,0);$ the dimension is $\frac{k}{2}+2$;
\item $\set{\varphi_{2p+1},1\ppq p\ppq \frac{k-4}{2};\varphi_{k-1};\theta_1;\theta_{-1}}$ if $k$ is even and $(c,d)=(0,0)$; the dimension is $\frac{k}{2}+1+\delta_{k,2}$.
\end{enumerate}
Moreover, when $k$ is odd or $k=2$, $\L^2(\HH^1(\Lambda))$ is spanned by:
\begin{enumerate}[(i)]
 \item $\set{\varphi_{2p+1},2\ppq p\ppq \frac{k-3}{2};\varphi_{k-1}}$ if $k$ is odd and $c=0$ and $d\neq 0$; the dimension is $\frac{k-3}{2}+\delta_{k,3}$;
 \item $\set{\varphi_{2p+1},2\ppq p\ppq \frac{k-3}{2}}$ if $k$ is odd and $(c,d)=(0,0)$; the dimension is $\frac{k-5}{2}+\delta_{k,3}$;
 \item $\set{\varphi_1,\theta_1,\theta_{-1}+d\theta_{-2}}$ if $k=2$ and  $c=0$ and $d\neq 0$; the dimension is $3$;
 \item $\set{\theta_1,\theta_{-1}}$ if $k=2$ and   $(c,d)=(0,0)$; the dimension is $2$.
  \end{enumerate}
\ep

As a consequence of Lemma \ref{lemma:q 2 when k=2} and Proposition \ref{prop:lc series q1}, we get the following result.

\bc\label{cor:classification q 2} Let $k\pgq 2$ be any integer and let $c$ and $d$ be non-zero elements in $K$. Then $Q(1\A)_1^k$,
$Q(1\A)_2^k(0,d)$, $Q(1\A)_2^k(c,d)$  are not stably equivalent of Morita type.
\ec

\br We still do not know whether $Q(1\A)_2^k(0,d)$ and $Q(1\A)_2^k(0,d')$ for $d\neq d'$ non-zero
are stably equivalent of Morita type 
or not or whether $Q(1\A)_2^k(c,d)$ and $Q(1\A)_2^k(c',d')$ for $\set{c,d}\neq \set{c',d'}$ with $cd\neq 0$ and $c'd'\neq 0$  are
stably equivalent of Morita type 
or not.

In fact, if $k$ is odd and $dd'\neq 0$, the Lie algebras $\HH^1(Q(1\A)_2^k(0,d))$ and  $\HH^1(Q(1\A)_2^k(0,d'))$ are isomorphic (in the remaining cases we do not know), and the isomorphism is given by 
\begin{align*}
&\varphi_t\mapsto \varphi_t'\text{ for $1\ppq t\ppq k-2$},&& \varphi_{k-1}\mapsto \frac{d'}{d}\varphi_{k-1}',\\&\psi\mapsto \frac{d}{d'}\psi',&& \theta_{-2}\mapsto \frac{d'}{d}\theta_{-2}'
\end{align*} with the obvious notations for the bases of the two Lie algebras.
\er

\end{document}